\documentclass[reqno,11pt]{article}
\usepackage{graphicx, amsmath, amsthm}
\usepackage{amssymb}

\input{stochbif.sty}

\begin{document}

%%%%%%%%%%%%%%%%%%%%%%%%%%%%%%%%%%%%%%%%%%%%%%%%%%%%%%%%%%%%%%%%%%%%%%%%%%%

\title{Metastability in Interacting Nonlinear \\
Stochastic Differential Equations~I: \\
From Weak Coupling to Synchronisation}
\author{Nils Berglund, Bastien Fernandez and Barbara Gentz}
\date{}

\maketitle

\begin{abstract}
\noindent
We consider the dynamics of a periodic chain of $N$ coupled overdamped
particles under the influence of noise. Each particle is subjected to a
bistable local potential, to a linear coupling with its nearest neighbours,
and to an independent source of white noise. The system shows a metastable
behaviour, which is characterised by the location and stability of its
equilibrium points. We show that as the coupling strength increases, the
number of equilibrium points decreases from $3^N$ to $3$. While for weak
coupling, the system behaves like an Ising model with spin-flip dynamics,
for strong coupling (of the order $N^2$), it synchronises, in the sense
that all particles assume almost the same position in their respective
local potential most of the time. We derive the exponential asymptotics for
the transition times, and describe the most probable transition paths
between synchronised states, in particular for coupling intensities below
the synchronisation threshold. Our techniques involve a centre-manifold
analysis of the desynchronisation bifurcation, with a precise control of
the stability of bifurcating solutions, allowing us to give a detailed
description of the system's potential landscape. 
\end{abstract}

\leftline{\small{\it Date.\/} November 21, 2006. Revised version, July 5, 2007.}
\leftline{\small 2000 {\it Mathematical Subject Classification.\/} 
37H20, 37L60 (primary), 37G40, 60K35 (secondary)}

%%%%%%%%%%%%%%%%%%%%%%%%%%%%%%%%%%%%%%%%%%%%%%%%%%%%%%%%%%%%%%%%%%%%%%%%%%%

\section{Introduction}
\label{sec_in}

%%%%%%%%%%%%%%%%%%%%%%%%%%%%%%%%%%%%%%%%%%%%%%%%%%%%%%%%%%%%%%%%%%%%%%%%%%%

Lattices of interacting deterministic multistable systems display a wide
range of interesting behaviours, due to the competition between local
dynamics and coupling between different sites. While for weak coupling,
they often exhibit spatial chaos (independent dynamics at the different
sites), for strong coupling they tend to display an organised collective
behaviour, such as synchronisation (see, for
instance~\cite{BelykhMosekilde96,ChowMalletVleck96,Johnston97,
NekorkinMakarovKazantsevVelarde97},
and~\cite{PRK,CML2004} for reviews).  

An important problem is to understand the effect of noise on such systems.
Noise can be used to model the effect of unresolved degrees of freedom, for
instance the influence of external heat reservoirs (see,
e.g.,~\cite{FKM,SL77,EPR,RBT00,RBT02}), which can induce currents through
the chain. The long-time behaviour of the system is described by its
invariant measure (assuming such a measure exists); however, for weak
noise, the dynamics often displays metastability, meaning that the
relaxation time towards the invariant measure is extremely long. 

Metastability has been studied extensively for particle systems with
stochastic dynamics. In these models, the transition from one metastable
state to another usually involves the gradual creation of a critical
droplet through small fluctuations, followed by a rapid transition to the
new state. The distributions of transition times, as well as the shapes of
critical droplets, have been investigated in detail (see in
particular~\cite{denHollander04,OlivieriVares05} for reviews, and
references therein).

In these lattice models, the local variables take only a finite number of
discrete values, which are independent of the interaction with other sites.
In the present paper, we consider by contrast a model with continuous
on-site variables. This leads to a system of interacting stochastic
differential equations (also called interacting diffusions, see, for
instance,~\cite{DawsonGaertner88} for a review of asymptotic properties in
the mean-field case). It turns out that while this system has a similar
behaviour to stochastic lattice models for weak coupling, the dynamics is
totally different for strong coupling: There are only $3$ equilibrium
configurations left, while the activation energy becomes extensive in the
number $N$ of particles. For large $N$, the system's behaviour is closer to
the behaviour of a Ginzburg--Landau partial differential equation with
noise (see, e.g.~\cite{EckmannHairer01,Rougemont02}).  The transition
between the strong-coupling and the weak-coupling regimes involves, as we
shall see, a sequence of symmetry-breaking bifurcations. Such bifurcations
have been studied, for instance, in~\cite{QinChen2004} for the
weak-coupling regime, and 
in~\cite{McNeil99,McNeil02,Watanabe93a,Watanabe93b} for systems of coupled
phase oscillators.

Our major aim is to determine the dependence of the transition times
between meta\-stable states, as well as the critical configurations, on the
coupling strength, on the whole range from weak to strong coupling. This
analysis requires a precise knowledge of the system's \lq\lq potential
landscape\rq\rq, in particular the number and location of its local minima
and saddles of index~$1$~\cite{FW,Sugiura96a,Kolokoltsov00,BEGK,BGK}. In
order to obtain this information, we will exploit the symmetry properties
of the system, using similar techniques as the ones developed in the
context of phase oscillators
in~\cite{AshwinSwift92,DionneGolubitskyStewart96a,
DionneGolubitskyStewart96b},
for instance. Our study also involves a centre-manifold analysis of the
desynchronisation bifurcation, which goes beyond existing results on
similar bifurcations because a precise control of the stability
of the bifurcating stationary points is required. 

This paper is organised as follows. Section~\ref{sec_res} contains the
precise description of the model and the statement of all results. After
introducing the model and describing its behaviour for weak and strong
coupling, in Section~\ref{ssec_symmetry} we examine the effect of
symmetries on the bifurcation diagram. A few special cases with small
particle number $N$ are illustrated in Section~\ref{ssec_smallN}. The
desynchronisation bifurcation for general $N$ is discussed in
Section~\ref{ssec_desybif}, while Section~\ref{ssec_otherbif} considers
further bifurcations of the origin. Finally, Section~\ref{ssec_stoch}
presents the consequences of these results for the stochastic dynamics of
the system. 

The subsequent sections contain the proofs of our results. The proof of
synchronisation at strong coupling is presented in Section~\ref{sec_prf},
while Section~\ref{sec_pf} introduces Fourier variables, which are used to
prove the results for $N=2$ and $N=3$, and for the centre-manifold analysis
of the desynchronisation bifurcation for general $N$.
Appendix~\ref{app_horse} gives a brief description of the analysis of the
weak-coupling regime, which uses standard techniques from symbolic
dynamics,
and Appendix~\ref{app_N4} contains a short description of the analysis of
the case $N=4$. 

The follow-up work~\cite{BFG06b} analyses in more detail the behaviour for
large particle number $N$. In that regime, we are able to control the
number
of stationary points in a much larger domain of coupling intensities,
including values far from the synchronisation threshold.

\subsection*{Acknowledgements}

Financial support by the French Ministry of Research, by way of the {\it
Action Concert\'ee Incitative (ACI) Jeunes Chercheurs, Mod\'elisation
stochastique de syst\`emes hors \'equilibre\/}, is gratefully acknowledged.
NB and BF thank the Weierstrass Institute for Applied Analysis and
Stochastics (WIAS), Berlin, for financial support and hospitality. BG
thanks the ESF Programme {\it Phase Transitions and Fluctuation Phenomena
for Random Dynamics in Spatially Extended Systems (RDSES)\/} for financial
support, and the Centre de Physique Th\'eorique (CPT), Marseille, for kind
hospitality. 

%%%%%%%%%%%%%%%%%%%%%%%%%%%%%%%%%%%%%%%%%%%%%%%%%%%%%%%%%%%%%%%%%%%%%%%%%%%

\section{Model and Results}
\label{sec_res}

%%%%%%%%%%%%%%%%%%%%%%%%%%%%%%%%%%%%%%%%%%%%%%%%%%%%%%%%%%%%%%%%%%%%%%%%%%%

\subsection{Definition of the Model}
\label{ssec_mod}

In our study of the influence of noise on lattice dynamical systems with
continuous on-site variables, we shall focus on a model which can serve as
a paradigm for such systems. It is based on coupled bistable ordinary
differential equations governed by the competition between individual
effects and spatial interactions. To name a few examples, it applies for
instance to chains of particles placed in a quartic potential and coupled
by springs~\cite{FM96}, to stimulus conduction in the myocardial
tissue~\cite{Keener87} and to certain chemical reactions~\cite{EN93}.
Another motivation for choosing this specific model is that it allows for a
quantitative description of phenomena, while being generic. In particular,
it quantifies changes from intensive to extensive properties when the
interaction strength increases.

The system is defined by the following ingredients:
\begin{itemiz}
\item	The periodic one-dimensional lattice is given by $\lattice=\Z/N\Z$,
where $N\geqs2$ is the number of particles.

\item	To each site $i\in\lattice$, we attach a real variable $x_i\in\R$,
describing the position of the $i$th particle. The configuration space is
thus $\cX=\R^\lattice$.

\item	Each particle feels a local bistable potential, given by 
\begin{equation}
\label{mod1}
U(\xi) = \frac14 \xi^4 - \frac12 \xi^2\;,
\qquad
\xi\in\R\;.
\end{equation}
The local dynamics thus tends to push the particle towards one of the two
stable positions $\xi=1$ or $\xi=-1$. The reason for this choice is that
$U(\xi)$ is the simplest possible double-well potential, which will be
responsible for metastability. Nevertheless, we expect this potential to
yield a behaviour which is, to a certain extent, generic among models in
its symmetry class. 

\item	Neighbouring particles in $\Lambda$ are coupled via a
discretised-Laplacian interaction, of intensity $\gamma/2$. Such a
nearest-neighbour coupling can be expected to yield very different
dynamics than mean-field models, for instance. 

\item	Each site is coupled to an independent source of noise, of
intensity $\sigma$. The sources of noise are described by independent
Brownian motions $\set{B_i(t)}_{t\geqs0}$.
\end{itemiz}

The system is thus described by the following set of coupled stochastic
differential equations, defining a diffusion on $\cX$:
\begin{equation}
\label{mod2}
\6x^\sigma_i(t) = f(x^\sigma_i(t))\6t 
+ \frac\cng2 \bigbrak{x^\sigma_{i+1}(t)-2x^\sigma_i(t)+x^\sigma_{i-1}(t)}
\6t
+ \sigma \6B_i(t)\;,
\qquad i\in\Lambda\;,
\end{equation}
where the local nonlinear drift is given by 
\begin{equation}
\label{mod3}
f(\xi) = -\nabla U(\xi) = \xi - \xi^3\;.
\end{equation}
For $\sigma=0$, the system~\eqref{mod2} is a gradient system of the form
$\dot x = -\nabla V_\gamma(x)$, with potential 
\begin{equation}
\label{mod5}
V_{\gamma}(x) = 
\sum_{i\in\lattice} U(x_i) + \frac\cng4 \sum_{i\in\lattice}
(x_{i+1}-x_i)^2\;.
\end{equation}

Note that the local potential $U(\xi)$ is invariant under the
transformation $\xi\mapsto-\xi$, implying that the local dynamics has no
preference for either positive or negative $\xi$. An interesting question
is how the results are affected by adding a symmetry-breaking term to
$U(\xi)$. This question will be the subject of future research. Some
preliminary studies indicate that several results, such as the presence of
synchronisation for strong coupling, the structure of the desynchronisation
bifurcation, and the qualitative behaviour for weak coupling, are not much
affected by the asymmetry. However, many details of the bifurcation
diagrams, as well as the metastable timescales and transition paths, will 
of course be quite different. 

%%%%%%%%%%%%%%%%%%%%%%%%%%%%%%%%%%%%%%%%%%%%%%%%%%%%%%%%%%%%%%%%%%%%%%%%%%%

\subsection{Potential Landscape and Metastability}
\label{ssec_pot}

The dynamics of the stochastic system depends essentially on the \lq\lq
potential landscape\rq\rq\ $V_\gamma$. 
More precisely, let 
\begin{equation}
\label{mod7}
\cS = \cS(\gamma)
= \setsuch{x\in\cX}{\nabla V_{\gamma}(x)=0}
\end{equation}
denote the set of stationary points of the potential. A point $x\in\cS$ is
said to be of type $(n_-,n_0,n_+)$ if the Hessian matrix of $V_\gamma$ at
$x$ has
$n_-$ negative, $n_+$ positive and $n_0=N-n_--n_+$ vanishing eigenvalues
(counting multiplicity). For each $k=0,\dots,N$, let 
$\cS_k=\cS_k(\gamma)$ denote the set of stationary points $x\in\cS$ which
are
of type $(k,0,N-k)$. For $k\geqs1$, these points are called\/
\emph{saddles of index $k$}\/, or simply\/ \emph{$k$-saddles}\/, while
$\cS_0$ is the set of strict local minima of $V_\gamma$. 

The stochastic system~\eqref{mod2} admits an invariant probability measure
with density proportional to $\e^{-2V_\gamma(x)/\sigma^2}$,  implying that
asymptotically, the system spends most of the time near the deepest minima
of $V_\gamma$. However, the invariant measure does not contain any
information on
the dynamics between these minima, nor on the way the equilibrium
distribution is approached.   Loosely speaking, for small noise intensity
$\sigma$ the stochastic system behaves in the following way~\cite{FW}: 
\begin{itemiz}
\item	A sample path $\set{x^\sigma(t)}_t$, starting in a point $x_0$
belonging to the deterministic basin of attraction $\cA(x^\star)$ of a
stationary point $x^\star\in\cS_0$, will first reach a small neighbourhood
of $x^\star$, in a time close to the time it would take a deterministic
solution to do so.  
\item	During an exponentially long time span, $x^\sigma(t)$ remains in
$\cA(x^\star)$, spending most of that time near $x^\star$, but making
occasional excursions away from the stationary point. 
\item	Sooner or later, $x^\sigma(t)$ makes a transition to (the
neighbourhood
of)  another stationary point $y^\star\in\cS_0$. During this transition,
the
sample path is likely to pass close to a saddle $s\in\cS_1$, whose unstable
manifolds converge to $x^\star$ and $y^\star$. In fact, the whole sample
path during the transition is likely to remain close to these unstable
manifolds. 
\item	After a successful transition, the sample path again spends an
  exponentially long time span in the basin of $y^\star$, until a similar
  transition brings it to another point of $\cS_0$ (which may or may not be
  different from $x^\star$).  
\end{itemiz}
If we ignore the excursions of the sample paths inside domains of
attraction, and consider only the transitions between local minima of the
potential, the stochastic process resembles a Markovian jump process on
$\cS_0$, with exponentially small transition probabilities, the only
relevant transitions being those between potential minima connected by a
$1$-saddle. 

Understanding the dynamics for small noise thus essentially requires
knowing the graph $\cG=(\cS_0,\cE)$, in which two vertices $x^\star,
y^\star\in\cS_0$ are connected by an edge $e\in\cE$ if and only if there is
a $1$-saddle $s\in\cS_1$ whose unstable manifolds converge to $x^\star$ and
$y^\star$. The mean transition time from $x^\star$ to $y^\star$ is of order
$\e^{2H/\sigma^2}$, where $H$ is the potential difference between $x^\star$
and the lowest saddle leading from $x^\star$ to $y^\star$.  

In our case, the potential $V_\gamma(x)$ being a polynomial of degree $4$
in $N$
variables, the set $\cS$ of stationary points admits at most $3^N$
elements. On the other hand, it is easy to see that $\cS$ always contains
at least the three points 
\begin{equation}
\label{pot1}
O = (0,\dots,0)\;, \qquad I^{\pm} = \pm(1,\dots,1)\;.
\end{equation}
Depending on the value of $\gamma$, the origin $O$ can be an $N$-saddle, or
a $k$-saddle for any odd $k$.  The points $I^{\pm}$ always belong to
$\cS_0$, in fact  we have for all $\gamma>0$
\begin{equation}
\label{pot2}
V_\gamma(x) > V_\gamma(I^+) = V_\gamma(I^-) = -\frac N4 \quad \forall
x\in\cX
\setminus\set{I^-,I^+}\;.
\end{equation}
Thus $I^-$ and $I^+$ are the most stable configurations of the system, and
also the local minima between which transitions take the longest time.

Among the many questions we can ask for the stochastic system, we shall
concentrate on the following:
\begin{itemiz}
\item	How long does the system typically take to make a transition from
$I^-$ to $I^+$, and how does the transition time depend on coupling
strength $\gamma$ and noise intensity $\sigma$? 
\item	How do typical paths for such a transition look like? 
\end{itemiz}

%%%%%%%%%%%%%%%%%%%%%%%%%%%%%%%%%%%%%%%%%%%%%%%%%%%%%%%%%%%%%%%%%%%%%%%%%%%

\subsection{Weak-Coupling Regime}
\label{ssec_small}

In the uncoupled case $\gamma=0$, we simply have 
\begin{equation}
\label{small1}
\cS(0) = \set{-1,0,1}^\lattice\;, 
\qquad
\abs{\cS(0)} = 3^N\;.
\end{equation}
Furthermore, $\cS_k(0)$ is the set of stationary points having exactly $k$
coordinates equal to $0$ (thus $\abs{\cS_k(0)} = \binom{N}{k}2^{N-k}$). In
particular, $\cS_0(0) = \set{-1,1}^\lattice$ has cardinality $2^N$. 

Hence the graph $\cG$ is an $N$-dimensional hypercube: Two vertices
$x^\star$, $y^\star\in\cS_0(0)$ are connected if and only if they differ in
exactly one component. Note that $V_0(x^\star)=-N/4$ for all local minima
$x^\star$, and $V_0(s)=-(N-1)/4$ for all $1$-saddles $s$, implying that all
nearest-neighbour transitions of the uncoupled system take the same time
(of order $\e^{1/2\sigma^2}$) on average.  

For small positive coupling intensity $0<\gamma\ll1$, the implicit-function
theorem guarantees that all stationary points depend analytically on
$\gamma$, without changing their type. In addition, the following result is
a direct consequence of standard results on H\'enon-like mappings of the
plane: 

\begin{prop}
\label{prop_sl1}
For any $N$, there exists a critical coupling $\gamma^\star(N)$ such that
for all $0\leqs\gamma<\gamma^\star(N)$, the stationary points
$x^\star(\gamma)\in\cS(\gamma)$ depend continuously on $\gamma$, without
changing their type. The critical coupling $\gamma^\star(N)$ satisfies 
\begin{equation}
\label{small2}
\inf_{N\geqs2}\gamma^\star(N) > \frac14\;.
\end{equation}
\end{prop} 

The proof is briefly discussed in Appendix~\ref{app_horse}, where we also 
provide slightly better lower bounds on $\inf_{N\geqs2}\gamma^\star(N)$. 
We
expect the critical coupling to be quite close to $1/4$, however. In
particular, we will show below that $\gamma^\star(2)=1/3$,
$\gamma^\star(3)=0.2701\dots$, and $\gamma^\star(4)=0.2684\dots$ 

\begin{figure}
\centerline{\includegraphics*[clip=true,height=60mm]{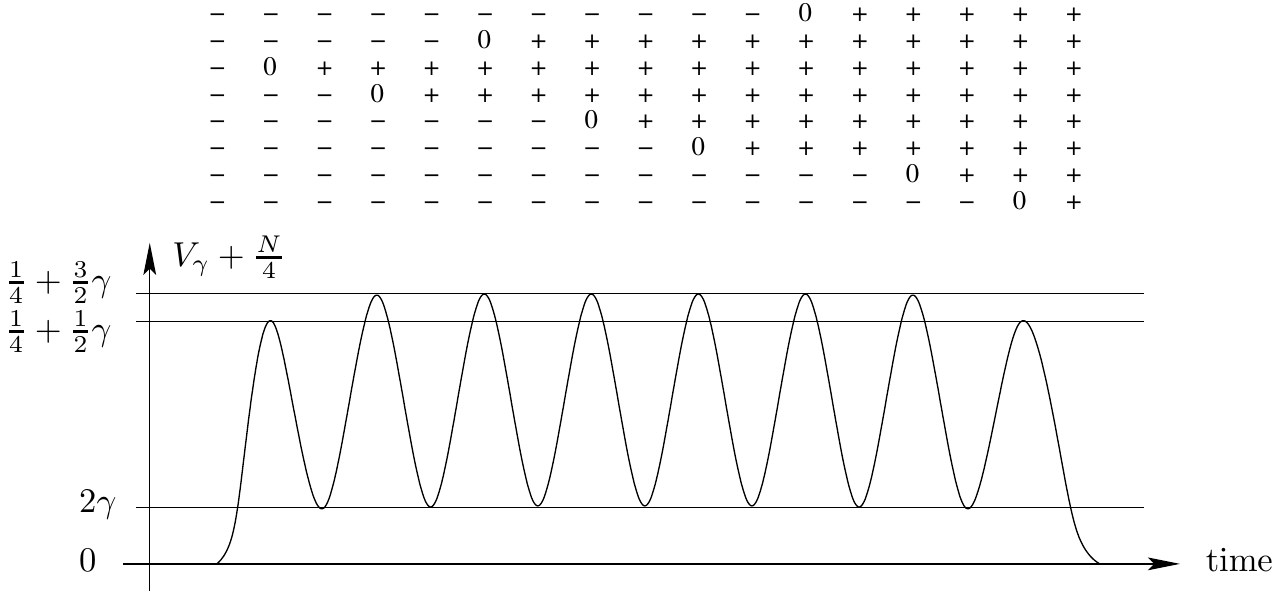}}
 \figtext{ }
 \caption[]
 {Example of an optimal transition path from $I^-$ to $I^+$ for weak
coupling.
 The upper half of the figure shows the local minima and $1$-saddles
visited
 during the transition, displayed vertically. For instance, the second
 column means that the first $1$-saddle visited is a perturbation of order
 $\gamma$ of the stationary point $(-1,-1,0,-1,-1,-1,-1,-1)$ present in
 absence of coupling. The lower half of the figure shows the value of the
 potential seen along the transition path.}
\label{fig_smallcoupling}
\end{figure}

%\newpage
Since any stationary point 
$x^\star(\gamma)=(x^\star_1(\gamma),\dots,x^\star_N(\gamma))\in\cS(\gamma)$
satisfies  $x^\star(\gamma)=x^\star(0)+\Order{\gamma}$, where each
$x^\star_i(0)$ is a stationary point of the local potential $U$, one has 
\begin{align}
\nonumber
V_\gamma(x^\star(\gamma)) 
&= V_\gamma(x^\star(0)) + \Order{\gamma^2} \\
&= V_0(x^\star(0)) + \frac\gamma4 \sum_{i=1}^N 
(x^\star_{i+1}(0)-x^\star_i(0))^2 + \Order{\gamma^2}\;. 
\label{small3}
\end{align}
To first order in $\gamma$, the potential's increase due to coupling
depends on the number of interfaces in the unperturbed configuration
$x^\star(0)\in\cS(0)$ (recall from~\eqref{small1} that the components of
$x^\star(0)$ only take values in $\set{-1,0,1}$). The dynamics of the
stochastic system is thus essentially the one of an Ising-spin system with
Glauber dynamics. Starting from the configuration $I^-$, the system reaches
with equal probability any configuration with one positive and $N-1$
negative spins. Then, however, it is less expensive to switch a spin
neighbouring the positive one than to switch a far-away spin, which would
create more interfaces. Thus the optimal transition from $I^-$ to $I^+$
consists in the growth of a \lq\lq droplet of $+$ in a sea of $-$ \rq\rq\
(\figref{fig_smallcoupling}).
To first order in the coupling intensity, all visited $1$-saddles except
the first and last one have the same potential value $-N/4 + 1/4 +
(3/2)\gamma + \Order{\gamma^2}$.  The energy required for the transition
is thus 
\begin{equation}
\label{small4}
V_\gamma\bigpar{(-1,-1,\dots,-1,0,1,1,\dots,1)+\Order{\gamma}} -
V_\gamma\bigpar{I^-}
= \frac14 + \frac32\gamma + \Order{\gamma^2}\;,
\end{equation} 
which is independent of the system size. Note that the situation is
different for lattices of dimension larger than $1$, in which the energy
increases with the surface of the droplet (cf.~\cite{denHollander04}). 

%%%%%%%%%%%%%%%%%%%%%%%%%%%%%%%%%%%%%%%%%%%%%%%%%%%%%%%%%%%%%%%%%%%%%%%%%%%

%\newpage
\subsection{Synchronisation Regime}
\label{ssec_sync}

For strong coupling $\gamma$, the situation is completely different than
for
weak coupling. Indeed, we have the following result (see
Section~\ref{sec_prf} for the proof):

\begin{prop}
\label{prop_sync1}
Let 
\begin{equation}
\label{synchro1}
\gamma_1 = \gamma_1(N) = \frac1{1-\cos(2\pi/N)}
= \frac{N^2}{2\pi^2}\biggbrak{1-\biggOrder{\frac1{N^2}}}\;.
\end{equation}
Then $\cS(\gamma) = \set{O,I^+,I^-}$ if and only if $\gamma\geqs\gamma_1$. 
Moreover, the origin $O$ is a $1$-saddle if and only if $\gamma>\gamma_1$,
and in that case its unstable manifold is contained in the diagonal 
\begin{equation}
\label{synchro2}
\cD = \setsuch{x\in\cX}{x_1=x_2=\dots=x_N}\;.
\end{equation}
\end{prop}

\begin{figure}
\centerline{\includegraphics*[clip=true,width=95mm]{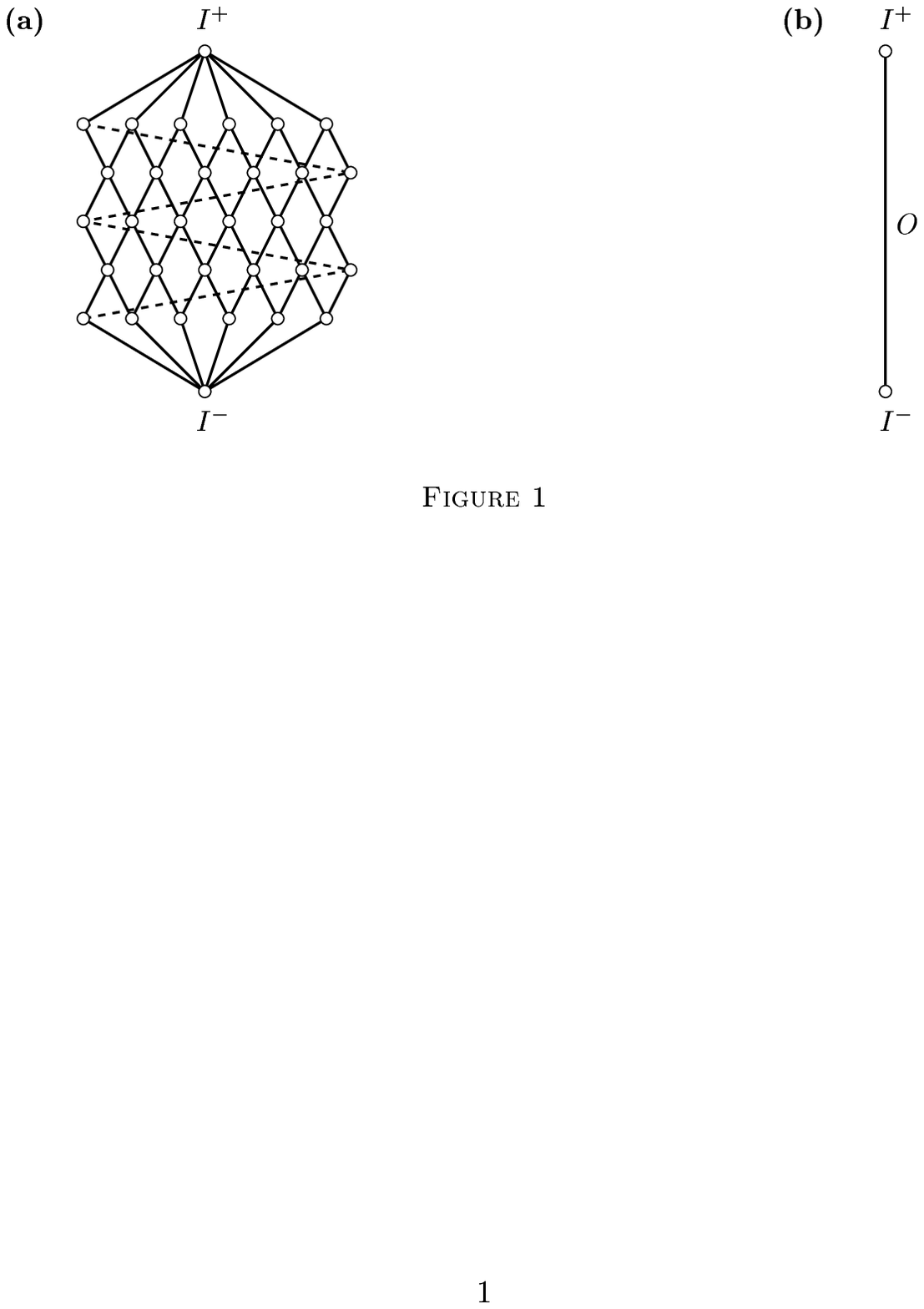}}
 \figtext{ }
 \caption[]
 {{\bf (a)}~Structure of the graph $\cG$ in the small-coupling regime
 $0<\gamma<\gamma^\star(N)$. Only edges belonging to optimal paths are
 shown: In the Ising-model analogy, these paths correspond to the flip of
 neighbouring spins. {\bf (b)}~The graph $\cG$ in the synchronisation
 regime $\gamma > \gamma_1(N)$.}
\label{fig_graphsext}
\end{figure}

We say that the deterministic system is\/ {\em synchronised}\/ for
$\gamma>\gamma_1$, in the sense that the diagonal is approached as
$t\to\infty$ for any initial state, meaning that all coordinates $x_i$ are
asymptotically equal. In other words, the system behaves as if all
particles coagulate in order to form one large particle of mass $N$. The
graph $\cG$ contains only two vertices $I^\pm$, connected by a single edge
(\figref{fig_graphsext}b). By extension, in the stochastic case we will say
that the system is synchronised whenever all  coordinates $x^\sigma_i$
remain close to each other most of the time with high probability.
Then transitions between $I^-$ and $I^+$ occur in a small neighbourhood of
the diagonal. 

In this case, the energy required for the transition is 
\begin{equation}
\label{sync2}
V_\gamma\bigpar{O} - V_\gamma\bigpar{I^-} = \frac N4\;,
\end{equation} 
which is extensive in the system size. Transitions now take a time of order
$\e^{N/2\sigma^2}$, and are thus much less frequent in the synchronisation
regime than in the low-coupling regime. 

It is remarkable that as the coupling $\gamma$ grows from $0$ to
$\gamma_1$, the number of stationary points decreases from $3^N$ to $3$.
The main purpose of this work is to elucidate in which way this transition
occurs, and how it affects the transition paths and times. 

%%%%%%%%%%%%%%%%%%%%%%%%%%%%%%%%%%%%%%%%%%%%%%%%%%%%%%%%%%%%%%%%%%%%%%%%%%%

%\newpage
\subsection{Symmetry Groups}
\label{ssec_symmetry}

The deterministic system $\dot x = -\nabla V_\gamma(x)$ is equivariant
(that is, $\nabla V_\gamma(g x)=g \nabla V_\gamma(x)$ for all $x\in\cX$) 
under three different types of symmetries $g$:
\begin{itemiz}
\item	Cyclic permutations (corresponding to rotations around the
diagonal), generated by
\begin{equation}
\label{sN2}
R(x_1,\dots,x_N) = (x_2,\dots,x_N,x_1)\;,
\end{equation}
as a consequence of the particles being identical.
\item	Reflection symmetries $R^kS$, $k=0,\dots,N-1$, where  
\begin{equation}
\label{sN3}
S(x_1,\dots,x_N) = (x_N,x_{N-1},\dots,x_1)\;,
\end{equation}
as a consequence of the interaction being isotropic.
\item	The inversion  
\begin{equation}
\label{sN4}
C(x_1,\dots,x_N) = -(x_1,\dots,x_N)\;,
\end{equation}
as a consequence of the local potential being even.
\end{itemiz}
The symmetries $R$ and $S$ generate the dihedral group $D_N$, which has
order $2N$,  for $N\geqs 3$, and the group $\Z_2$ for $N=2$. For
$N\geqs3$,  the symmetries $R$, $S$ and $C$ generate a group of order $4N$,
which we shall denote $G=G_N=D_N\times\Z_2$. For $N=2$ the symmetry group
is the Klein four-group $G_2=\Z_2\times\Z_2$, which as order $4$.

The set of stationary points $\cS(\gamma)$, as well as each set
$\cS_k(\gamma)$ of $k$-saddles, are invariant under $G$. Thus $G$ acts as a
group of transformations on $\cX$, on $\cS(\gamma)$, and on each
$\cS_k(\gamma)$. We will use a few concepts from elementary group theory:
\begin{itemiz}
\item	For $x\in\cX$, the\/ {\em orbit}\/ of $x$ is the set
$O_x=\setsuch{gx}{g\in G}$. 
\item	For $x\in\cX$, the\/ {\em isotropy group}\/ or {\em stabiliser}\/
of $x$ is the set
$C_x=\setsuch{g\in G}{gx=x}$.
\item	The\/ {\em fixed-point set}\/ of a subgroup $H$ of $G$ is the set
$\Fix(H)=\setsuch{x\in\cX}{hx=x \;\forall h\in H}$.
\end{itemiz}  
The following facts are well known:
\begin{itemiz}
\item	For any $x$, the isotropy group $C_x$ is a subgroup of $G$ and
$\abs{C_x}\abs{O_x}=\abs{G}$ (where $\abs{X}$ denotes the cardinality of
a finite set $X$).

\item	For any $g\in G$ and $x\in\cX$, we have $C_{gx} = gC_xg^{-1}$,
so that the isotropy groups of all the points of a given orbit are
conjugate. 

\item	For any subgroup $H$ of $G$ and any $g\in G$,
$\Fix(gHg^{-1})=g\Fix(H)$.  
\end{itemiz}

These facts allow us to limit the study to one point of each orbit, to one
subgroup in each conjugacy class, and to one type of conjugated fixed-point
set. 
For small $N$, this reduction often suffices to completely determine all
stationary points of the system, while for larger $N$, it at least helps to
classify the stationary points.

%%%%%%%%%%%%%%%%%%%%%%%%%%%%%%%%%%%%%%%%%%%%%%%%%%%%%%%%%%%%%%%%%%%%%%%%%%%

\subsection{Small Lattices}
\label{ssec_smallN}

\begin{figure}[t]
\centerline{\includegraphics*[clip=true,width=110mm]{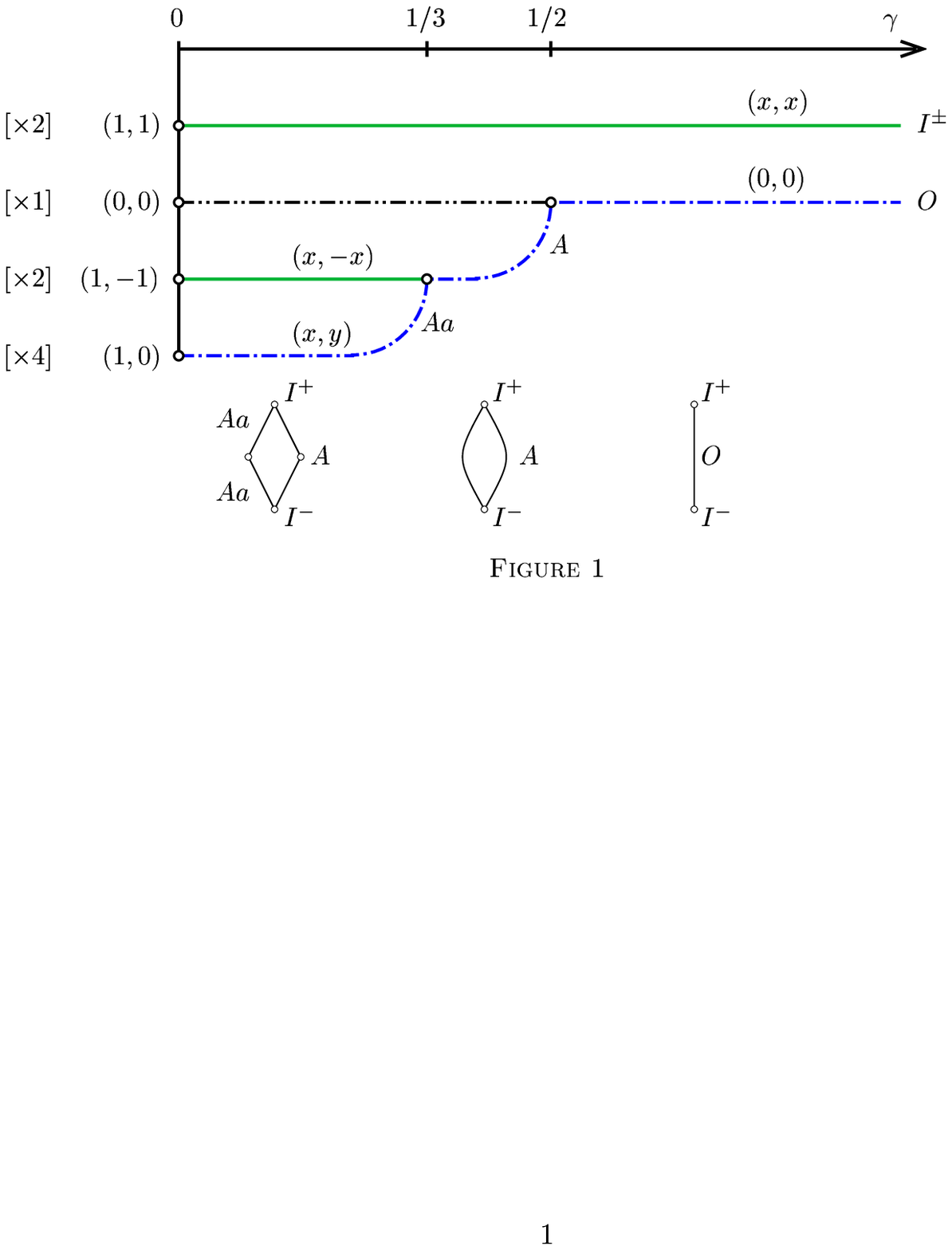}}
 \figtext{ }
 \caption[]
 {Bifurcation diagram for the case $N=2$ and associated graphs $\cG$. Only
 one stationary point is shown for each orbit of the symmetry group $G$.
 The cardinality of the orbit is shown in square brackets. 
 Full lines represent local minima of the potential, while dash--dotted
 lines with $k$ dots represent $k$-saddles.}
\label{fig_bif2}
\end{figure}

\begin{table}[t]
\begin{center}
\begin{tabular}{|l|l|l|l|}
\hline
\vrule height 12pt depth 6pt width 0pt
$z^\star$ & $O_{z^\star}$ & $C_{z^\star}$ & $\Fix(C_{z^\star})$ \\
\hline
\vrule height 12pt depth 6pt width 0pt
$(0,0)$ & $\set{(0,0)}$ & $G_2=\Z_2\times\Z_2$ & $\set{(0,0)}$ \\
\vrule height 6pt depth 6pt width 0pt
$(1,1)$ & $\set{(1,1),(-1,-1)}$ & $\Z_2=\set{\id,S}$ &
$\set{(x,x)}_{x\in\R}=\cD$ \\
\vrule height 6pt depth 6pt width 0pt
$(1,-1)$ & $\set{(1,-1),(-1,1)}$ & $\set{\id,CS}$ & $\set{(x,-x)}_{x\in\R}$
\\
\vrule height 6pt depth 8pt width 0pt
$(1,0)$ & $\set{\pm(1,0),\pm(0,1)}$ & $\set{\id}$ &
$\set{(x,y)}_{x,y\in\R}=\cX$ \\
\hline 
\end{tabular} 
\end{center}
\caption[]
{Stationary points $z^\star\in\cS(0)$, their orbits, their isotropy groups
and
the corresponding fixed-point sets in the case $N=2$.}
\label{table_N=2}
\end{table}

We now consider some particular cases for illustration. The following
applies to the three stationary points that are always present:
\begin{itemiz}
\item	For the origin $O$, $O_O = \set{O}$, $C_O = G$ and
$\Fix(C_O)=\set{O}$. 
\item	For the global minima $I^\pm$, $O_{I^+}=O_{I^-}=\set{I^-,I^+}$,
$C_{I^\pm}=D_N$ and $\Fix(C_{I^\pm})=\cD$, where $\cD$ is the diagonal
defined in~\eqref{synchro2}. 
\end{itemiz}

\subsubsection*{Case $N=2$}

For $N=2$, we have $R=S$ and the symmetry group is
$G_2=\set{\id,S,C,CS}$. In the uncoupled case, the set of stationary points
$\cS(0)$ is partitioned into four orbits, as shown in
Table~\ref{table_N=2}. 
\figref{fig_bif2} indicates how the stationary points evolve as the
coupling
increases (the proof is given in Proposition~\ref{prop_f2}). We see that
stationary points keep the same type of symmetry as the coupling
strength $\gamma$ increases, and sometimes merge with a stationary point of
higher symmetry. 

Below the bifurcation diagram, we show how the corresponding graphs $\cG$
change as the coupling intensity $\gamma$ varies. The
two-dimen\-sional hypercube (i.e., the square), present for weak coupling,
transforms into a graph with two vertices, connected by two edges, as the
$1$-saddles labelled $Aa$ undergo a pitchfork bifurcation at
$\gamma=1/3=\gamma^\star(2)$. For $1/3<\gamma<1/2$, the points with
$(x,-x)$-symmetry, labelled $A$, are $1$-saddles, representing the points
of maximal potential height on the two optimal transition paths from
$I^-$ to $I^+$. At $\gamma=1/2=\gamma_1(2)$, the $1$-saddles undergo
another
pitchfork bifurcation, this time with the origin $O$, which becomes the
only transition gate in the strong-coupling regime.

The value of the potential on the bifurcating branches is found to be 
\begin{align}
\nonumber
V_\gamma(A) &= -\frac12 (1-2\gamma)^2\;,\\
V_\gamma(Aa) &= -\frac12 (1-2\gamma)^2 + \frac14 (1-3\gamma)^2\;.
\label{smallN1}
\end{align} 
Thus in the case $\gamma<1/3=\gamma^\star(2)$, the typical time for a
transition from $I^-$ to $I^+$ will be of order
$\e^{(1+2\gamma+\gamma^2)/2\sigma^2}$, while
for $1/3<\gamma<1/2$, it is of order $\e^{4\gamma(1-\gamma)/\sigma^2}$,
and for $\gamma>1/2$ it is of order $\e^{1/\sigma^2}$. 

\begin{figure}[t]
\centerline{\includegraphics*[clip=true,width=110mm]{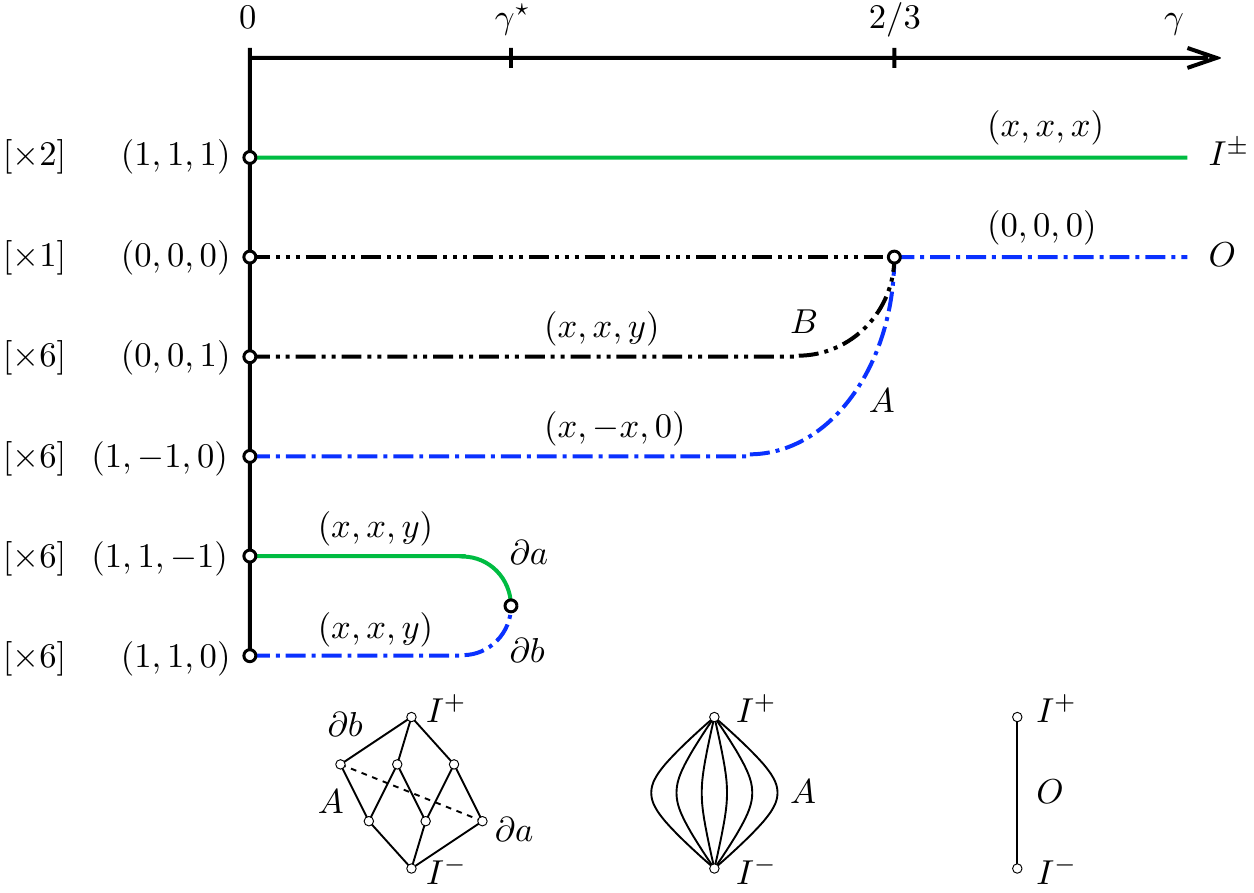}}
% \figtext{ }
 \caption[]
 {Bifurcation diagram for the case $N=3$ and associated graphs $\cG$. Only
 one stationary point is shown for each orbit of the symmetry group $G$. 
 The saddle--node bifurcation value is  $\gamma^\star = \gamma^\star(3) = 
 (\sqrt{3+2\sqrt3}-\sqrt3)/3 = 0.2701\dots$}
\label{fig_bif3}
\end{figure}

\begin{table}
\begin{center}
\begin{tabular}{|l|l|l|l|}
\hline
\vrule height 12pt depth 6pt width 0pt
$z^\star$ & $O_{z^\star}$ & $C_{z^\star}$ & $\Fix(C_{z^\star})$ \\
\hline
\vrule height 12pt depth 6pt width 0pt
$(0,0,0)$ & $\set{(0,0,0)}$ & $G_3$ & $\set{(0,0,0)}$ \\
\vrule height 6pt depth 6pt width 0pt
$(1,1,1)$ & $\set{(1,1,1),(-1,-1,-1)}$ & $D_3$ &
$\set{(x,x,x)}_{x\in\R}=\cD$ \\
\vrule height 6pt depth 6pt width 0pt
$(1,-1,0)$ & $\set{\pm(1,-1,0),\pm(-1,0,1),\pm(0,1,-1)}$ & $\set{\id,CRS}$
&
$\set{(x,-x,0)}_{x\in\R}$ \\
\vrule height 6pt depth 6pt width 0pt
$(0,0,1)$ & $\set{\pm(0,0,1),\pm(0,1,0),\pm(1,0,0)}$ & $\set{\id,RS}$ &
$\set{(x,x,y)}_{x,y\in\R}$ \\
\vrule height 6pt depth 6pt width 0pt
$(1,1,-1)$ & $\set{\pm(1,1,-1),\pm(1,-1,1),\pm(-1,1,1)}$ & $\set{\id,RS}$ &
$\set{(x,x,y)}_{x,y\in\R}$ \\
\vrule height 6pt depth 8pt width 0pt
$(1,1,0)$ & $\set{\pm(1,1,0),\pm(1,0,1),\pm(0,1,1)}$ & $\set{\id,RS}$ &
$\set{(x,x,y)}_{x,y\in\R}$ \\
\hline 
\end{tabular} 
\end{center}
\caption[]
{Stationary points $z^\star\in\cS(0)$, their orbits, their isotropy groups
and
the corresponding fixed-point sets in the case $N=3$.}
\label{table_N=3}
\end{table}

\subsubsection*{Case $N=3$}

For $N=3$, at zero coupling the set $\cS(0)$ of stationary points is
partitioned into six orbits, as shown in Table~\ref{table_N=3}. Their
evolution as the coupling increases is shown in \figref{fig_bif3} (for a
proof, see Proposition~\ref{prop_f3}). The new feature in this case is that
two orbits disappear in a saddle--node bifurcation at
$\gamma=\gamma^\star(3)=(\sqrt{3+2\sqrt3}-\sqrt3)/3=0.2701\dots$, instead
of merging with stationary points of higher symmetry. This accounts for the
rather drastic transformation of the graph $\cG$ from a $3$-cube for
$\gamma<\gamma^\star(3)$ to a graph with two vertices joined by $6$ edges
for $\gamma^\star(3)<\gamma<2/3=\gamma_1(3)$. The potential on the
$A$-branches has value
\begin{equation}
\label{smallN2}
V_\gamma(A) = -\frac12\Bigpar{1-\frac32\gamma}^2\;.
\end{equation}
Thus in the case $\gamma<2/3$, the typical time for a
transition from $I^-$ to $I^+$ will be of order
$\e^{(2+12\gamma-9\gamma^2)/4\sigma^2}$, while
for $\gamma>2/3$ it is of order $\e^{3/2\sigma^2}$. 

\begin{figure}[t]
\centerline{\includegraphics*[clip=true,width=120mm]{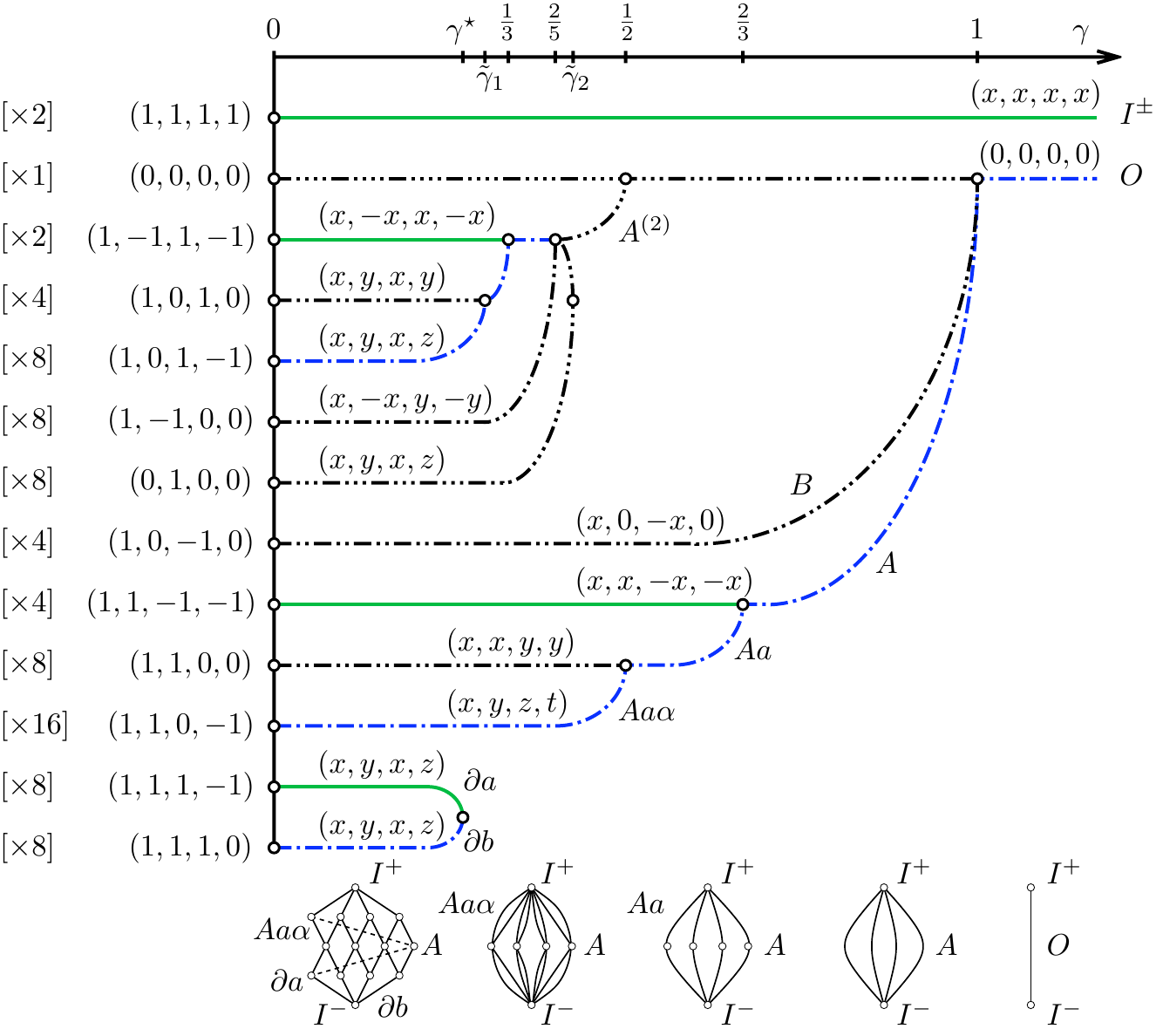}}
% \figtext{ }
 \caption[]
 {Bifurcation diagram for the case $N=4$ and associated graphs $\cG$ (only
 edges corresponding to optimal transition paths are shown). Only one
 stationary point is shown for each orbit of the symmetry group $G$.  The
 saddle--node bifurcation values are  $\gamma^\star = \gamma^\star(4) =
 0.2684\dots$ and $\tilde\gamma_2 = 0.4004\dots$, while $\tilde\gamma_1 =
 (3\sqrt2 - 2)/7 = 0.3203\dots$}
\label{fig_bif4}
\end{figure}

\subsubsection*{Case $N=4$}

\figref{fig_bif4} shows the results of a similar analysis for $N=4$.  The
determination of this bifurcation diagram, which relies partly on the
numerical study of the roots of some polynomials, is outlined in
Appendix~\ref{app_N4}.  As in the previous cases, a certain number of
stationary points emerge from the origin as the coupling intensity
decreases below the synchronisation threshold. In the present case, the
desynchronisation bifurcation occurs at $\gamma=\gamma_1(4)=1$, and there
are four $1$-saddles, labelled $A$, and four $2$-saddles, labelled $B$,
emerging from the origin. Two more symmetry-breaking bifurcations affect
the $A$-branches, finally resulting in $1$-saddles without any symmetry. A
second branch, labelled $A^{(2)}$, bifurcates from the origin at
$\gamma=1/2$. We do not show the corresponding stationary points in the
graphs, because they appear always to correspond to non-optimal transition
paths. 

The examples discussed here give some flavour of how the transition from
weak coupling to synchronisation occurs in the general case. In the sequel,
we will mainly describe the desynchronisation bifurcation occurring at
$\gamma=\gamma_1(N)$, which can be analysed for arbitrary~$N$. 

%%%%%%%%%%%%%%%%%%%%%%%%%%%%%%%%%%%%%%%%%%%%%%%%%%%%%%%%%%%%%%%%%%%%%%%%%%%

\subsection{Desynchronisation Bifurcation}
\label{ssec_desybif}

In this section, we consider the behaviour for general particle number
$N\geqs3$, 
and coupling intensities $\gamma$ slightly below the synchronisation
threshold 
$\gamma_1(N)$. Our aim is to obtain precise information on
\begin{itemiz}
\item 	the exact number and stability type of stationary points;
\item	the potential difference between global minima and saddles of index
$1$, 
which governs the transition time between these minima;
\item	the spatial shape of the critical configurations (i.e., the saddles
of index 
$1$), which are the highest energy configurations reached during a typical 
transition. 
\end{itemiz}

\subsubsection*{Even particle number}

We start by considering the case of even particle number $N$. The main
result is the 
following theorem, which will be proved in Section~\ref{sec_cm}. 

\begin{theorem}[Desynchronisation bifurcation, even particle number]
\label{thm_desync1}
Assume that $N$ is even. Then there exist $\delta=\delta(N)>0$ and points
$A=A(\gamma)$ and $B=B(\gamma)$ in $\cX$ such that for
$\gamma_1-\delta < \gamma < \gamma_1$, the set of stationary points
$\cS(\gamma)$ has cardinality $2N+3$, and can be decomposed as follows: 
\begin{align}
\nonumber
\cS_0 &= O_{I^+} = \set{I^+,I^-}\;, \\
\nonumber
\cS_1 &= O_{A} = \set{A,RA,\dots,R^{N-1}A}\;, \\
\nonumber
\cS_2 &= O_{B} = \set{B,RB,\dots,R^{N-1}B}\;, \\
\cS_3 &= O_{O} = \set{O}\;. 
\label{desybif1}
\end{align}
Furthermore, the value $V_\gamma(A)$ of the potential is the same on 
all $1$-saddles, and is determined by
\begin{equation}
\label{desybif4B}
\frac{V_\gamma(A)}N = 
\begin{cases}
\vrule height 12pt depth 16pt width 0pt
-\dfrac14(1-\gamma)^2 
& \text{if $N=4$\;,} \\
\vrule height 16pt depth 10pt width 0pt
-\dfrac16 \bigpar{1-\gamma/\gamma_1}^2 
+ \bigOrder{(1-\gamma/\gamma_1)^3}
& \text{if $N\geqs 6$\;,}
\end{cases}
\end{equation}
\end{theorem}

This result tells us that as the coupling intensity $\gamma$ decreases
below the synchronisation threshold $\gamma_1$, there are exactly $N$
saddles of index $1$ and $N$ saddles of index $2$ bifurcating from the
origin. All saddles of the same index belong to the same orbit of the
symmetry group $G$. In fact, if the location of one saddle is known,
the locations of all other saddles of the same type are obtained by
cyclic permutation of its coordinates. 

Relation~\eqref{desybif4B} shows that the height of the $1$-saddles
decreases as the coupling intensity $\gamma$ decreases. A simple
topological argument shows that the unstable manifolds of each $1$-saddle
converge to $I^+$ and $I^-$. Thus the transition time from $I^-$ to $I^+$
will be governed by the potential difference $V_\gamma(A)-V_\gamma(I^-)$. 

\begin{remark}
The proof of Theorem~\ref{thm_desync1} actually yields more precise
information on the coordinates of the saddles $A$ and $B$, which can be
summarised as follows. 

\begin{itemiz}

\item The saddles $A$ and $B$ satisfy special symmetries, as shown in
Table~\ref{table_desync} and~\figref{fig_poly}. Namely, their coordinates
admit a mirror symmetry, as well as a mirror symmetry with sign change, the
two symmetry axes being perpendicular. The details depend on whether
$N\in4\N$ or $N\in4\N+2$, because this affects the number of coordinates
which may lie on the symmetry axes. 

\begin{table}
\begin{center}
\begin{tabular}{|l|l|l|l|}
\hline
\vrule height 12pt depth 6pt width 0pt
$N$ & $x$ & $C_x$ & $\Fix(C_x)$ \\
\hline
\vrule height 12pt depth 6pt width 0pt
$4L$ & $A$ & $\Dbar_2$ &
$(x_1,\dots,x_L,x_L,\dots,x_1,-x_1,\dots,-x_L,-x_L,\dots,-x_1)$ \\
\vrule height 8pt depth 6pt width 0pt
     & $B$ & $\Dbar_2'$ &
$(x_1,\dots,x_L,\dots,x_1,0,-x_1,\dots,-x_L,\dots,-x_1,0)$ \\
\hline 
\vrule height 12pt depth 6pt width 0pt
$4L+2$ & $A$ & $\Dbar_2$ &
$(x_1,\dots,x_{L+1},\dots,x_1,-x_1,\dots,-x_{L+1},\dots,-x_1)$ \\
\vrule height 8pt depth 6pt width 0pt
     & $B$ & $\Dbar_2'$ &
$(x_1,\dots,x_L,x_L\dots,x_1,0,-x_1,\dots,-x_L,-x_L,\dots,-x_1,0)$ \\
\hline 
\vrule height 12pt depth 6pt width 0pt
$2L+1$ & $A$ & $\langle CRS\rangle$ &
$(x_1,\dots,x_L,-x_L,\dots,-x_1,0)$ \\
\vrule height 8pt depth 6pt width 0pt
     & $B$ & $\langle RS\rangle$ &
$(x_1,\dots,x_L,x_L,\dots,x_1,x_0)$ \\
\hline 
\end{tabular} 
\end{center}
\caption[]
{Symmetries of the stationary points bifurcating from the origin at
$\gamma=\gamma_1$. Here $L$ denotes an integer such that $N=4L$, $N=4L+2$,
or $N=2L+1$ depending on the value of $N \pmod{4}$. 
The isotropy groups for even $N$ are $\Dbar_2=\langle CS,
R^{N/2}S\rangle$ and $\Dbar_2'=\langle CRS, R^{N/2+1}S\rangle$, where
$\langle
g_1,\dots,g_m \rangle$ denotes the group generated by $\set{g_1, \dots,
g_m}$.}
\label{table_desync}
\end{table}

\begin{figure}
\centerline{\includegraphics*[clip=true,width=120mm]{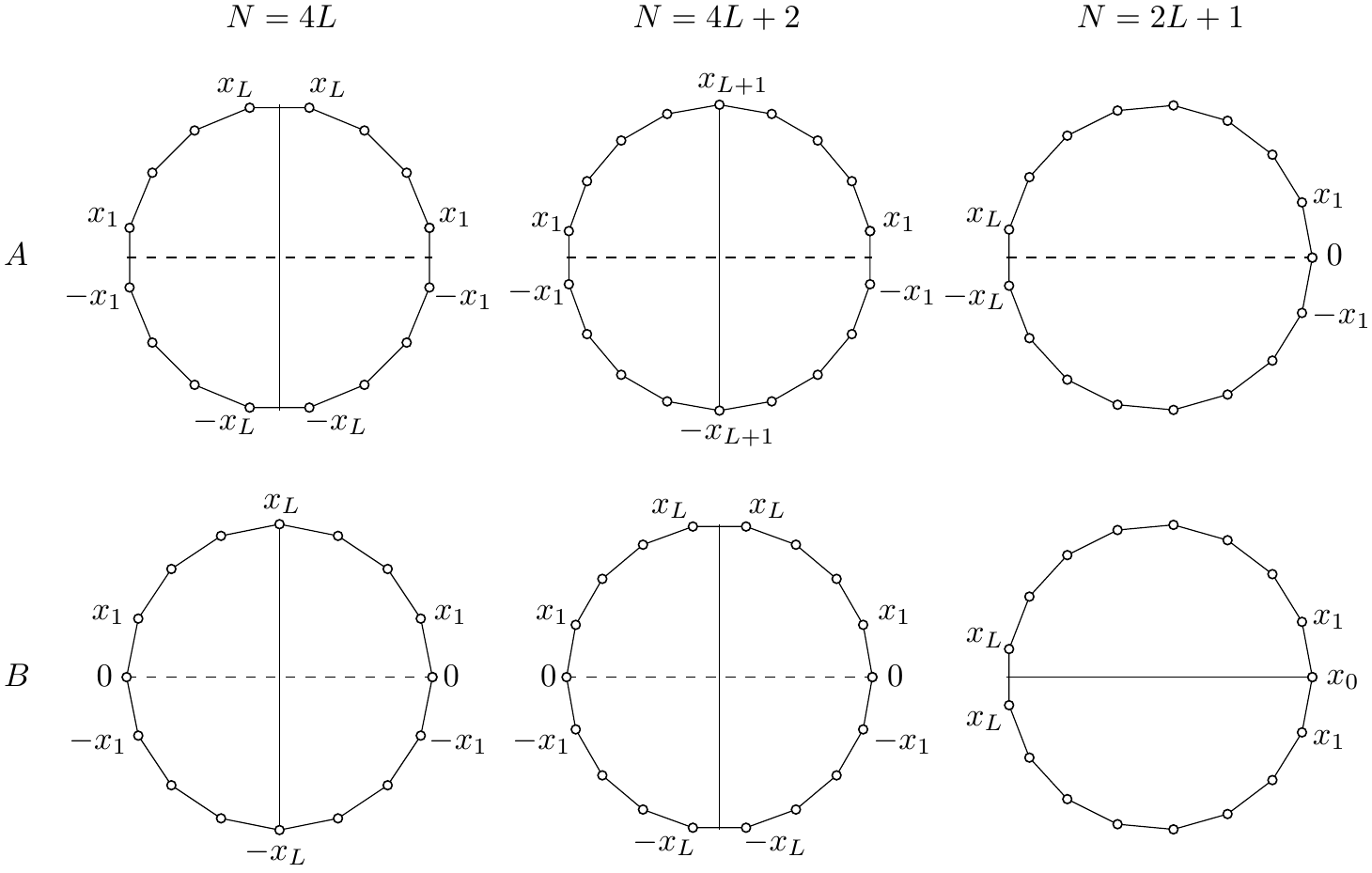}}
 \figtext{ }
 \caption[]
 {Symmetries of the stationary points $A$ and $B$ bifurcating from the
origin at
 $\gamma=\gamma_1$. Full lines represent symmetry axes, broken lines
 represent symmetry axes with sign change.}
\label{fig_poly}
\end{figure}

\item	The value $V_\gamma(B)$ of the potential on the $2$-saddles
satisfies 
\begin{equation}
\label{desybif4C}
0 < 
\frac{V_\gamma(B) -V_\gamma(A)}N \leqs 
\bigOrder{(1-\gamma/\gamma_1)^{N/2}}\;,
\end{equation}
indicating that the potential between the two types of saddles becomes
flatter as the particle number increases. 

\item	The components of $A=A(\gamma)$ and $B=B(\gamma)$ are given by  
\begin{align}
\nonumber
A_j(\gamma) &= \frac2{\sqrt3}\, \sqrt{1-\gamma/\gamma_1}
\sin \biggpar{\frac{2\pi}N \bigpar{j-\tfrac12}} 
+ \bigOrder{1-\gamma/\gamma_1}\;,\\
B_j(\gamma) &= \frac2{\sqrt3}\, \sqrt{1-\gamma/\gamma_1}
\sin \biggpar{\frac{2\pi}N j} 
+ \bigOrder{1-\gamma/\gamma_1}\;,
\label{desybif3}
\end{align}
(see \figref{fig_sym}), 
except in the case $N=4$, where  
\begin{equation}
\label{desybif2}
A(\gamma)=\sqrt{1-\gamma}\,(1,1,-1,-1) 
\qquad
\text{and} 
\qquad
B(\gamma)=\sqrt{1-\gamma}\,(1,0,-1,0)\;.
\end{equation}

\begin{figure}
\centerline{\includegraphics*[clip=true,width=130mm]{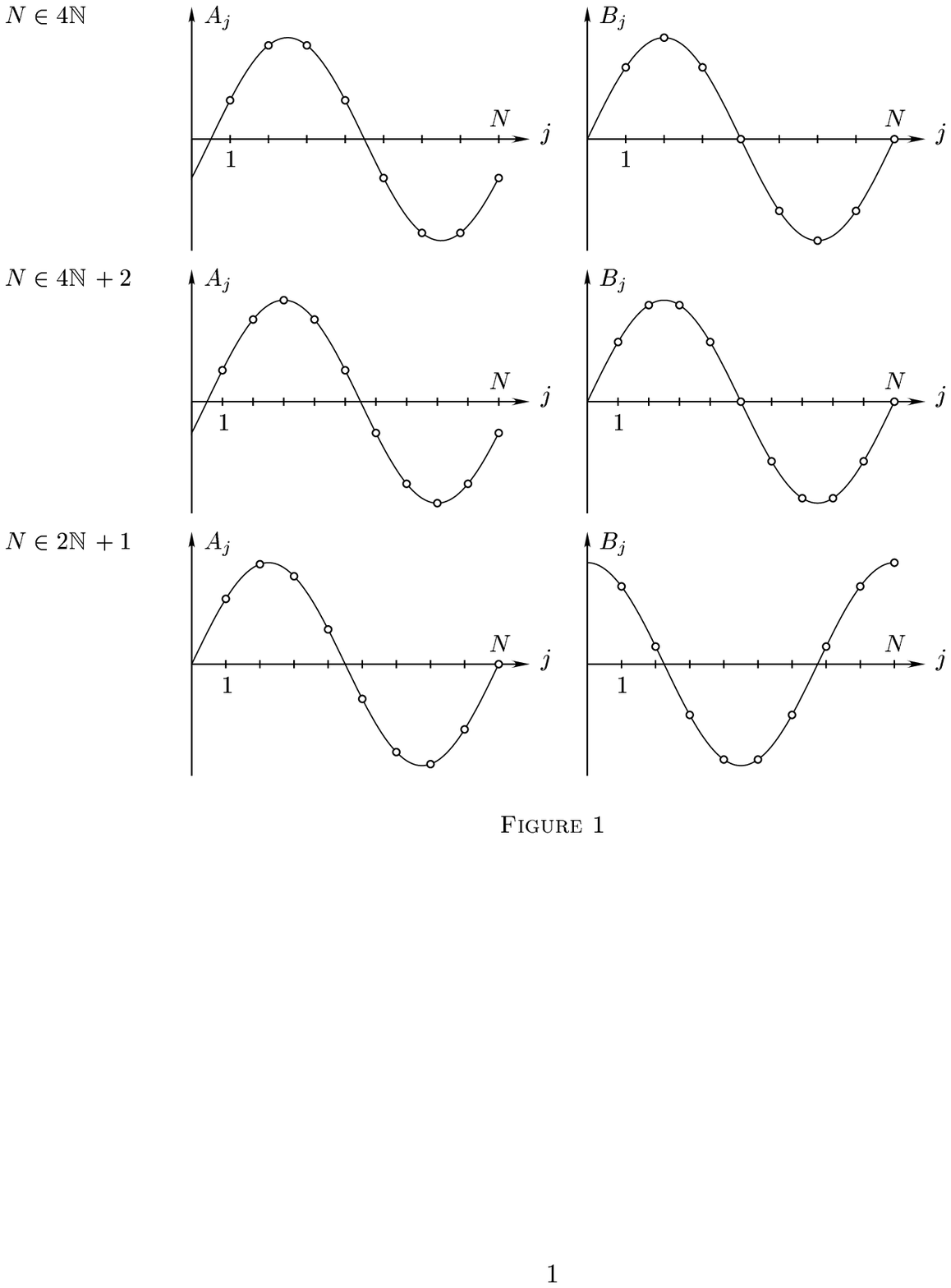}}
 \figtext{ }
 \caption[]
 {Components of the stationary points $A$ and $B$ bifurcating from the
 origin at $\gamma=\gamma_1$, shown for the three different cases $N=4L$,
 $N=4L+2$ and $N=2L+1$.}
\label{fig_sym}
\end{figure}

\item	For any $\gamma\in[0,\gamma_1)$, there exist  stationary
points $A(\gamma)$ and $B(\gamma)$, satisfying the symmetries indicated in
Table~\ref{table_desync}. The first $N/2$ components of $A(\gamma)$ and the
first $N/2-1$ components of $B(\gamma)$ are always positive, and 
\begin{align}
\nonumber
\lim_{\gamma\to0}A(\gamma) &= (1,1,\dots,1,1,-1,-1,\dots,-1,-1)\;,\\
\lim_{\gamma\to0}B(\gamma) &= (1,1,\dots,1,0,-1,-1,\dots,-1,0)\;.
\label{desybif4}
\end{align}
We do not claim that $A(\gamma)$ and $B(\gamma)$ are continuous in
$\gamma$ everywhere, though we expect them to be so. What we do prove is
that for any $\gamma$, there is at least one stationary point with the
appropriate symmetry (and coordinates of the indicated signs). We
also know that $A$ and $B$ depend continuously on $\gamma$ for $\gamma$
near $0$ and near $\gamma_1$. We cannot exclude, however, the presence of
saddle--node bifurcations in between.

\item	We know that the points $A(\gamma)$ are local minima near
$\gamma=0$. Thus these stationary points must undergo at least one
secondary
bifurcation as $\gamma$ decreases. For symmetry reasons, we expect that, as
in the case $N=4$ (see~\figref{fig_bif4}), there are two successive
symmetry-breaking bifurcations affecting the $1$-saddles: First, the mirror
symmetry with sign change is broken, then the remaining mirror symmetry
is broken as well. 

\item	 The error terms in~\eqref{desybif4B} and~\eqref{desybif3} may
depend on $N$. The technique we employ here does not allow for an optimal
control of the $N$-dependence of these error terms and of $\delta(N)$.
However, in~\cite{BFG06b}, we obtain such a control in the limit of large
$N$, using different techniques. 
\end{itemiz}
\end{remark}

\subsubsection*{Odd particle number}

In the case of odd particle number $N$, the results are slightly weaker,
since we are not able to obtain such a precise control on the number of
$1$-saddles and $2$-saddles created in the desynchronisation bifurcation. 
The main result is the following theorem, which will also be proved in
Section~\ref{sec_cm}. 

\goodbreak

\begin{theorem}[Desynchronisation bifurcation, odd particle number]
\label{thm_desybif2}
Assume that $N$ is odd. Then there exists $\delta=\delta(N)>0$ such that
for $\gamma_1-\delta < \gamma < \gamma_1$, the set of stationary points
$\cS(\gamma)$ has cardinality $4\ell N+3$, for some $\ell\geqs1$. All
stationary points are saddles of index $0, 1, 2$ or $3$, with 
\begin{align}
\nonumber
\cS_0 &= O_{I^+} = \set{I^+,I^-}\;, \\
\cS_3 &= O_{O} = \set{O}\;,
\label{desybif5}
\end{align}
and $\abs{\cS_1}=\abs{\cS_2}=2\ell N$. 
The value $V_\gamma(A)$ of the potential on any saddle $A$ of index $1$
satisfies 
\begin{equation}
\label{desybif5555}
\frac{V_\gamma(A)}N = 
-\dfrac16 \bigpar{1-\gamma/\gamma_1}^2 
+ \bigOrder{(1-\gamma/\gamma_1)^3}\;.
\end{equation}
\end{theorem}

In fact, we expect that $\ell=1$, so that there are exactly $2N$ saddles
of index $1$ and $2N$ saddles of index $2$. This would lead to the
following conjecture. 

\begin{conjecture}[Desynchronisation bifurcation, odd particle number]
\label{con_desync2}
For odd $N$ and  $\gamma_1-\delta < \gamma < \gamma_1$, the set of
stationary points $\cS(\gamma)$ has cardinality $4N+3$. There exist
points $A=A(\gamma)$ and $B=B(\gamma)$ in $\cX$ such that $\cS$ can be
decomposed as follows: 
\begin{align}
\nonumber
\cS_0 &= O_{I^+} = \set{I^+,I^-}\;, \\
\nonumber
\cS_1 &= O_{A} = \set{A,RA,\dots,R^{N-1}A,-A,-RA,\dots,-R^{N-1}A}\;, \\
\nonumber
\cS_2 &= O_{B} = \set{B,RB,\dots,R^{N-1}B,-B,-RB,\dots,-R^{N-1}B}\;, \\
\cS_3 &= O_{O} = \set{O}\;. 
\label{desybif5B}
\end{align}
\end{conjecture}

This conjecture would follow as a consequence of a much simpler
conjecture on the behaviour of certain coefficients in a centre-manifold
expansion, which can be computed iteratively, see Section~\ref{sec_cm}. We
know that it is true for $N=3$, and it can be checked by direct computation
for the first few values of $N$. Numerically, we checked the validity of
the conjecture for $N$ up to $101$. In~\cite{BFG06b},
we show that the conjecture is also true for sufficiently large $N$.

\begin{remark}
Again, the proof of the theorem also yields more precise information on
the location of the saddles. 
\begin{itemiz}
\item	The potential difference between the $2$-saddles and
$1$-saddles satisfies 	
\begin{align}
0 < \frac{V_\gamma(B) -V_\gamma(A)}N 
\leqs \bigOrder{(1-\gamma/\gamma_1)^N}\;.
\label{desybif777}
\end{align}

\item 	If the conjecture is true, then the components of $A=A(\gamma)$ and
$B=B(\gamma)$ are of the form
\begin{align}
\nonumber
A_j(\gamma) &= \frac2{\sqrt3}\, \sqrt{1-\gamma/\gamma_1}
\sin \biggpar{\frac{2\pi}N j} 
+ \bigOrder{1-\gamma/\gamma_1}\;,\\
B_j(\gamma) &= \frac2{\sqrt3}\, \sqrt{1-\gamma/\gamma_1} 
\cos \biggpar{\frac{2\pi}N j} 
+ \bigOrder{1-\gamma/\gamma_1}\;,
\label{desybif6B}
\end{align}
and satisfy the symmetries indicated in
Table~\ref{table_desync}. 

\item	Independently of the validity of the conjecture, for any
$\gamma\in[0,\gamma_1)$, there exist stationary points $A(\gamma)$ and
$B(\gamma)$, satisfying the symmetries indicated in
Table~\ref{table_desync}. The first $(N-1)/2$ components of $A(\gamma)$ are
strictly positive, and  
\begin{equation}
\lim_{\gamma\to0}A(\gamma) = (1,\dots,1,-1,\dots,-1,0)\;.
\label{desybif7B}
\end{equation}
\end{itemiz}
\end{remark}

%%%%%%%%%%%%%%%%%%%%%%%%%%%%%%%%%%%%%%%%%%%%%%%%%%%%%%%%%%%%%%%%%%%%%%%%%%%

\subsection{Subsequent Bifurcations}
\label{ssec_otherbif}

We shall show in Section~\ref{sec_prf} that the origin undergoes further
bifurcations at 
\begin{equation}
\label{dsybif8}
\gamma = \gamma_M = \frac1{1-\cos(2\pi M/N)}\;,
\qquad
2\leqs M\leqs \intpart{N/2}\;,
\end{equation}
in which the index of the origin $O$ increases by $2$ (except for the case
where $N$ is even and $M=N/2$, where the index increases by $1$), and new
saddles $A^{(M)}$ and $B^{(M)}$ of index $2M-1$ and $2M$ are created. 
Consequently these saddles are not important for the stochastic dynamics,
and we shall not provide a detailed analysis here. We briefly mention a few
properties of these bifurcations, which we will prove in~\cite{BFG06b}
to hold for sufficiently large $N/M$:
\begin{itemiz}
\item	The number of newly created stationary points is given by
$4N/\gcd(N,2M)$,
where $\gcd(N,2M)$ denotes the greatest common divisor of $N$ and $2M$.
\item	The number of sign changes of $x_j$ as a function of $j$ for these
new stationary points is equal to $2M$. $M$ can therefore be considered as
a\/
\emph{winding number}\/.
\item	If $N$ and $M$ are coprime, the new stationary points $x$ satisfy
the
symmetries shown in Table~\ref{table_desync}, while for other $M$ they
belong to larger isotropy subgroups. 
\end{itemiz}

\begin{example}
\label{ex_desybif}
For $N=8$, the origin bifurcates four times as $\gamma$ decreases from
$+\infty$ to $0$. We show the
symmetries of the bifurcating stationary points in
Table~\ref{table_exdesync}. They are obtained in the following way:
\begin{itemiz}
\item	Compute the eigenvectors of the Hessian of the potential at the
origin;
\item	Determine the corresponding isotropy subgroups of $G_8$;
\item	Write the equation $\dot z=-\nabla V_\gamma(z)$ restricted to the
fixed-point set of each isotropy subgroup, and study the bifurcations of
the
origin in each restricted system. 
\end{itemiz}

\begin{table}
\begin{center}
\begin{tabular}{|c|c|c|c|c|c|}
\hline
\vrule height 14pt depth 6pt width 0pt
$M$ & $\gamma_M$ & $\gcd(N,2M)$ & $A^{(M)}$ & $B^{(M)}$ \\
\hline
\vrule height 14pt depth 6pt width 0pt
$1$ & $2+\sqrt2$ & $2$ &  
$(x,y,y,x,-x,-y,-y,-x)$ & $(x,y,x,0,-x,-y,-x,0)$ \\
\vrule height 12pt depth 6pt width 0pt
$2$ & $1$ & $4$ & 
$(x,x,-x,-x,x,x,-x,-x)$ & $(x,0,-x,0,x,0,-x,0)$ \\
\vrule height 12pt depth 6pt width 0pt
$3$ & $2-\sqrt2$ & $2$ & 
$(x,-y,-y,x,-x,y,y,-x)$ & $(x,-y,x,0,-x,y,-x,0)$ \\
\vrule height 12pt depth 6pt width 0pt
$4$ & $1/2$ & $8$ & 
$(x,-x,x,-x,x,-x,x,-x)$ &  \\
\hline 
\end{tabular} 
\end{center}
\caption[]
{Fixed-point sets of the stationary points bifurcating from the origin
for $N=8$, for different winding numbers $M$. Points of winding number
$M=1$ and $M=3$ actually have the same fixed point spaces, but we choose
the signs of the components in such a way that $x$ and $y$ always have the
same sign for the actual stationary points.}
\label{table_exdesync}
\end{table}

For instance, for winding number $M=1$ and orbits of type $A$, we obtain 
\begin{equation}
\label{desybif9}
\begin{split}
\dot x &= (1-\tfrac32\gamma) x + \tfrac12\gamma y - x^3\;,\\
\dot y &= \tfrac12\gamma x + (1-\tfrac12\gamma) y - y^3\;. 
\end{split}
\end{equation} 
The origin bifurcates for $\gamma=2\pm\sqrt2$. An analysis of the
linearised system shows that for $\gamma$ slightly smaller than $2+\sqrt2$,
the new stationary points must lie in the quadrants
$\setsuch{(x,y)}{x>0,y>0}$ and $\setsuch{(x,y)}{x<0,y<0}$. These quadrants,
however, are invariant under the flow of $\dot x=-\nabla V_\gamma(x)$,
since, e.g., $\dot x>0$ if $x=0$ and $y>0$, and $\dot y>0$ if $y=0$ and
$x>0$. Hence the points created in the bifurcation remain in these
quadrants. The points created at $\gamma=2-\sqrt2$, which correspond to the
winding number $M=3$, lie in the complementary quadrants
$\setsuch{(x,y)}{x>0,y<0}$ and $\setsuch{(x,y)}{x>0,y<0}$. 

For $\gamma=0$, the only points with the appropriate symmetry are 
\begin{align}
\nonumber
A^{(1)}(0) &= (1,1,1,1,-1,-1,-1,-1)\;,
&
B^{(1)}(0) &= (1,1,1,0,-1,-1,-1,0)\;, \\
\nonumber
A^{(2)}(0) &= (1,1,-1,-1,1,1,-1,-1)\;,
&
B^{(2)}(0) &= (1,0,-1,0,1,0,-1,0)\;, \\
\nonumber
A^{(3)}(0) &= (1,-1,-1,1,-1,1,1,-1)\;,
&
B^{(3)}(0) &= (1,-1,1,0,-1,1,-1,0)\;, \\
A^{(4)}(0) &= (1,-1,1,-1,1,-1,1,-1)\;. 
\label{largeN9b}
\end{align} 
Note that the cases $M=2$ and $M=4$ are obtained by concatenation of
multiple copies of stationary points existing for $N=4$ and $N=2$,
respectively. 
\end{example}

%%%%%%%%%%%%%%%%%%%%%%%%%%%%%%%%%%%%%%%%%%%%%%%%%%%%%%%%%%%%%%%%%%%%%%%%%%%

\subsection{Stochastic Case}
\label{ssec_stoch}

We return now to the behaviour of the system of stochastic differential
equations
\begin{equation}
\label{stoch1}
\6x^\sigma_i(t) = f(x^\sigma_i(t))\6t 
+ \frac\cng2 \bigbrak{x^\sigma_{i+1}(t)-2x^\sigma_i(t)+x^\sigma_{i-1}(t)}
\6t
+ \sigma \6B_i(t)\;, \qquad i\in\Lambda\;.
\end{equation}
Recall that our main goal is to characterise the noise-induced transition
from the configuration $I^-=(-1,-1,\dots,-1)$ to the configuration
$I^+=(1,1,\dots,1)$. In particular, we are interested in the time needed
for this transition to occur, and in the shape of the critical 
configuration, i.e., the configuration of highest energy reached in the
course of the transition. 

Since the probability of a stochastic process in continuous space hitting a
given point is typically zero, we have to work with small neighbourhoods of
the relevant configurations.
Given a Borel set $\cA\subset\cX$, and an initial condition
$x_0\in\cX\setminus\cA$, we denote by $\tauhit(\cA)$ the\/
\emph{first-hitting time of $\cA$}\/ 
\begin{equation}
\label{stoch2}
\tauhit(\cA) = \inf\setsuch{t>0}{x^\sigma(t)\in\cA}\;.
\end{equation}
Similarly, for an initial condition $x_0\in\cA$, we denote by
$\tauexit(\cA)$ the\/ \emph{first-exit time from $\cA$}\/ 
\begin{equation}
\label{stoch3}
\tauexit(\cA) = \inf\setsuch{t>0}{x^\sigma(t)\notin\cA}\;.
\end{equation}
We can now formulate our main results, which are similar in spirit to
Theorems~3.2.1 and~4.2.1 of~\cite{denHollander04}. 

\begin{theorem}[Stochastic case, synchronisation regime]
\label{thm_stoch1}
Assume that the coupling strength satisfies
$\gamma>\gamma_1=(1-\cos(2\pi/N))^{-1}$. We fix radii $0<r<R<1/2$, and
denote by $\tau_+=\tauhit(\cB(I_+,r))$ the first-hitting time of a ball
$\cB(I^+,r)$ of radius $r$ around $I^+$. Then for any initial condition
$x_0\in\cB(I^-,r)$, any $N\geqs2$ and any $\delta>0$,
\begin{equation}
\label{stoch4}
\lim_{\sigma\to0}
\bigprobin{x_0}{\e^{(N/2-\delta)/\sigma^2} < \tau_+ <
\e^{(N/2+\delta)/\sigma^2}} = 1
\end{equation}
and 
\begin{equation}
\label{stoch4b}
\lim_{\sigma\to0} \sigma^2 \log 
\expecin{x_0}{\tau_+} = \frac N2\;.
\end{equation}
Furthermore, let $\tau_O = \tauhit(\cB(O,r))$, and let 
\begin{equation}
\label{stoch5}
\tau_- = \inf\setsuch{t>\tauexit(\cB(I^-,R))}{x_t\in\cB(I^-,r)}
\end{equation}
be the time of first return to the small ball $\cB(I^-,r)$ around
$I^-$ after leaving the larger ball $\cB(I^-,R)$. Then 
\begin{equation}
\label{stoch6}
\lim_{\sigma\to0}
\bigpcondin{x_0}{\tau_O < \tau_+}{\tau_+ < \tau_-} = 1\;.
\end{equation}
\end{theorem}

\begin{figure}
\centerline{\includegraphics*[clip=true,height=40mm]{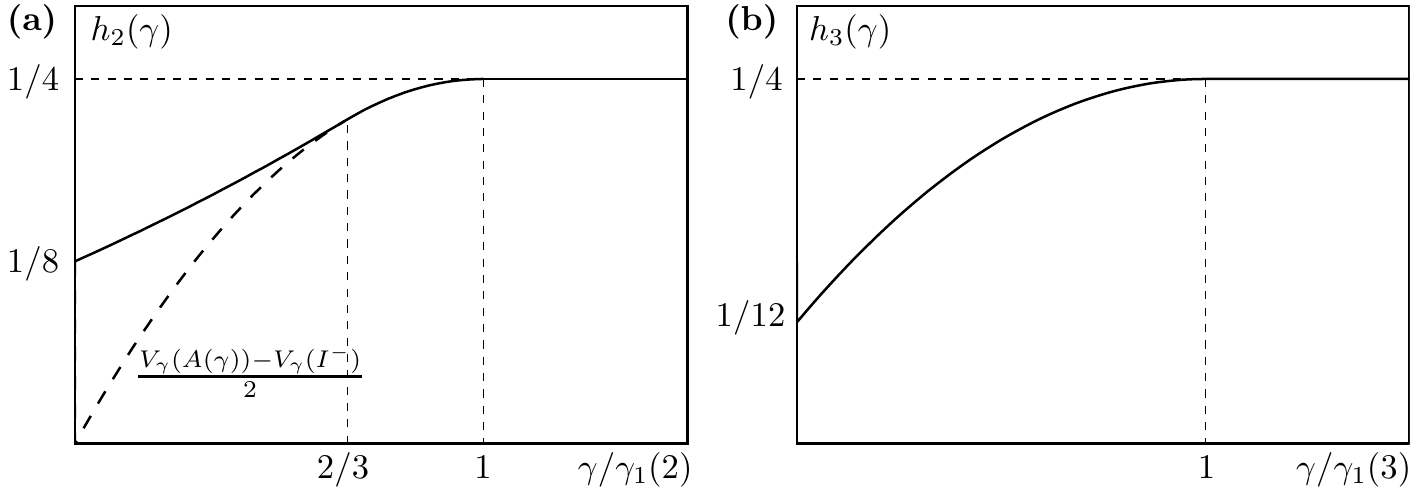}}
 \figtext{ }
 \caption[]
 {The normalised potential differences $h_N(\gamma)$ for $N=2$ and $N=3$.
 The broken curve for $N=2$ shows the potential difference
 $(V_\gamma(A)-V_\gamma(I^-))/2$, which is smaller than $h_N(\gamma)$ for
 $\gamma/\gamma_1<2/3$ (compare~\eqref{smallN1}), because for these
 parameter values, the $1$-saddle is the point labelled $Aa$.  For general
 particle number $N$, we know that $h_N(\gamma)$ behaves like
 $1/4-c_N(1-\gamma/\gamma_1)^2$ as $\gamma\nearrow\gamma_1$, and that
 $h_N(0)=1/4N$.}
\label{fig_HN}
\end{figure}

We will omit the proof of this result, which is a direct consequence of 
Proposition~\ref{prop_sync1} and standard results in Wentzell--Freidlin
theory, see for instance \cite{FW,Sugiura96a,Kifer}. 

Relations~\eqref{stoch4} and~\eqref{stoch4b} mean that the transition from
$I^-$ to $I^+$ typically takes a time of order $\e^{N/2\sigma^2}$. 
Relation~\eqref{stoch6} shows that, provided a transition from $\cB(I^-,r)$
to $\cB(I^+,r)$ is observed, the process is likely to pass close to the
saddle at the origin on its way from one potential well to the other one.
The origin thus plays the r\^ole of the critical configuration in the
synchronisation regime. 

\begin{theorem}[Stochastic case, desynchronised regime]
\label{thm_stoch2}
Assume $\gamma<\gamma_1$. 
Let $r, R$ and $\tau_+$ and $\tau_-$ be defined as in
Theorem~\ref{thm_stoch1}, and fix an initial condition $x_0\in\cB(I^-,r)$.
Let $h_N(\gamma)=(V_\gamma(A)-V_\gamma(I^-))/N$ denote the normalised
potential difference between $1$-saddles and potential minima at $I^-$. 
Then $h_N(\gamma)$ satisfies
\begin{equation}
\label{stoch7}
h_N(\gamma) = \frac14 - c_N \bigpar{1-\gamma/\gamma_1}^2 +
\bigOrder{(1-\gamma/\gamma_1)^3}
\qquad
\text{as $\gamma\nearrow\gamma_1$\;,}
\end{equation}
where $c_2=c_4=1/4$ and $c_N=1/6$ for $N=3$ and all $N\geqs5$, 
and one has 
\begin{equation}
\label{stoch8}
\lim_{\sigma\to0}
\bigprobin{x_0}{\e^{(2Nh_N(\gamma)-\delta)/\sigma^2} < \tau_+ <
\e^{(2Nh_N(\gamma)+\delta)/\sigma^2}} = 1\;, 
\end{equation}
and
\begin{equation}
\label{stoch8b}
\lim_{\sigma\to0} \sigma^2 \log 
\expecin{x_0}{\tau_+} = 2Nh_N(\gamma)\;.
\end{equation}
Furthermore, assume that either $N$ is even, or $N$ is odd and
Conjecture~\ref{con_desync2} holds. Let  
\begin{equation}
\label{stoch9}
\tau_A = \tauhit\Bigpar{\bigcup_{g\in G}\cB(gA(\gamma),r)}\;. 
\end{equation}
Then there exists a $\delta=\delta(N)>0$ such that for
$\gamma_1-\delta < \gamma < \gamma_1$
\begin{equation}
\label{stoch10}
\lim_{\sigma\to0}
\bigpcondin{x_0}{\tau_A < \tau_+}{\tau_+ < \tau_-} = 1\;.
\end{equation}
\end{theorem}

The proof is again standard, and we will omit it here. 

Relations~\eqref{stoch8} and~\eqref{stoch8b} mean that the transition from
$I^-$ to $I^+$ typically takes a time of order
$\e^{2Nh_N(\gamma)/\sigma^2}$, while relation~\eqref{stoch10} shows that
the set of critical configurations is the orbit of $A(\gamma)$. 

We conclude the statement of results with a few comments:
\begin{itemiz}
\item	The very precise results in~\cite{BEGK,BGK} allow in principle for
a more precise control of the expected transition time $\tau_+$ than the
exponential asymptotics given in~\eqref{stoch4b} and~\eqref{stoch8b}.
However, these results cannot be applied directly to our case, because they
assume some non-degeneracy condition to hold for the potential, which
excludes symmetries as in the system we are investigating. Therefore, in
the present work we content ourselves with the exponential asymptotics. A
finer analysis is feasible, and we plan to provide such a further study.

\item	We have seen that the potential difference between the $1$-saddles
$A$ and the $2$-saddles $B$ becomes smaller and smaller as the particle
number increases. As a consequence, for any given noise intensity
$\sigma>0$, the transition paths become less localised as the particle
number $N$ increases. This reflects the fact that the system becomes
translation-invariant in the large-$N$ limit. 

\item	As the coupling intensity $\gamma$ decreases, the configurations of
highest energy reached in the course of a transition become more and more
inhomogeneous along the chain, and we expect them to converge to
configurations of the form $(1,\dots,1,0,-1,\dots,-1)$ in the uncoupled
limit. In addition, several local minima and saddles should appear by
saddle--node bifurcations along the optimal transition path, thereby
increasing the number of metastable states of the system. 
\end{itemiz}

%%%%%%%%%%%%%%%%%%%%%%%%%%%%%%%%%%%%%%%%%%%%%%%%%%%%%%%%%%%%%%%%%%%%%%%%%%%

\newpage
\section{Lyapunov Functions and Synchronisation}
\label{sec_prf}

%%%%%%%%%%%%%%%%%%%%%%%%%%%%%%%%%%%%%%%%%%%%%%%%%%%%%%%%%%%%%%%%%%%%%%%%%%%

We now turn to the proofs of the statements made in Section~\ref{sec_res},
and  start by introducing a few notations.  The \lq\lq interaction
part\rq\rq\ of the potential is proportional to  
\begin{equation}
\label{pp8}
W(x) = \frac12 \sum_{i\in\Lambda} (x_i-x_{i+1})^2 
= \frac12 \norm{x-Rx}^2 
= \pscal{x}{\Sigma x}\;,
\end{equation}
where $\Sigma$ is the symmetric matrix  $\Sigma = \one -
\frac12(R+\transpose{R})$. Hence, the potential $V_\gamma(x)$ can be
written as 
\begin{equation}
\label{ppA1}
V_\gamma(x) = -\frac12 \pscal{x}{(\one-\gamma\Sigma)x} 
+ \frac14 \sum_{i\in\lattice} x_i^4\;.
\end{equation}
The eigenvectors of $\Sigma$ are of the form
$v_k=\transpose{(1,\omega^k,\dots,\omega^{(N-1)k})}$, $k=0,\dots,N-1$,
where $\omega=\e^{2\pi\icx/N}$, with eigenvalues  $1-\cos(2\pi k/N)$. This
implies in particular that the Hessian of the potential at the origin,
which is given by $\gamma\Sigma-\one$, has eigenvalues $-\lambda_k$, where 
\begin{equation}
\label{ppA2}
\lambda_k = \lambda_{-k} = 1 - \gamma \biggpar{1-\cos\frac{2\pi k}N}
= 1 - \frac\gamma{\gamma_k}\;.
\end{equation} 
The origin is a $1$-saddle for $\gamma>\gamma_1$. As $\gamma$ decreases,
the
index of the origin increases by $2$ each time $\gamma$ crosses one of the
$\gamma_k$, until it becomes an $N$-saddle at
$\gamma=\gamma_{\intpart{N/2}}$. 

We now show that for $\gamma>\gamma_1$ (i.e., $\lambda_1<0$), $W(x)$ is a
Lyapunov function for the deterministic system  $\dot x=-\nabla
V_\gamma(x)$. 

\begin{prop}
\label{prop_ps1}
For any initial condition $x_0$, the solution $x(t)$ of  $\dot x=-\nabla
V_\gamma(x)$ satisfies 
\begin{equation}
\label{ppA3}
\dtot{}{t} W(x(t)) \leqs 2 \bigpar{1-\gamma/\gamma_1} W(x(t)) 
- \frac1N W(x(t))^2\;.
\end{equation}
As a consequence, if $\gamma$ is strictly larger than $\gamma_1$, then
$x(t)$  converges exponentially fast to the diagonal, and thus the only
equilibrium points of the system are $O$ and $I^\pm$. 
\end{prop}
\begin{proof}
We first observe that the relation 
\begin{equation}
\label{ps1:1}
f(x_i) - f(x_{i+1}) = (x_i - x_{i+1})
\bigbrak{1-(x_i^2+x_ix_{i+1}+x_{i+1}^2)}
\end{equation}
allows us to write 
\begin{equation}
\label{ps1:2}
\dtot{}{t} (x-Rx) = \Pi(x,\gamma) (x-Rx)\;,
\end{equation}
where  $\Pi(x,\gamma) = \one - \gamma \Sigma - D(x)$. Here $D(x)$ is a
diagonal
matrix, whose $i$th entry is given by $x_i^2+x_ix_{i+1}+x_{i+1}^2$, and can
be bounded below by $\frac14(x_i-x_{i+1})^2$. 
It follows 
\begin{align}
\nonumber
\dtot{}t W(x(t)) 
&= \pscal{x-Rx}{\dtot{}t (x-Rx)} \\
\nonumber
&= \pscal{x-Rx}{\Pi(x,\gamma) (x-Rx)} \\
\nonumber
&\leqs \pscal{x-Rx}{(\one - \gamma \Sigma) (x-Rx)} 
- \frac14 \sum_{i\in\lattice} (x_i-x_{i+1})^4 \\
&\leqs \lambda_1 \norm{x-Rx}^2 - \frac1{4N} \norm{x-Rx}^4
\label{ps1:3}
\end{align}
by Cauchy-Schwartz. This implies~\eqref{ppA3}. If $\gamma\geqs\gamma_1$,
then $W(x(t))$ converges to zero as $t\to\infty$ for all initial
conditions, which implies that all stationary points $x^\star$ must satisfy
$W(x^\star)=0$, and thus lie on the diagonal. The only stationary points on
the diagonal, however, are $O$ and $I^\pm$. 
\end{proof}

\begin{proof}[{\sc Proof of Proposition~\ref{prop_sync1}}]
The assertion on the unstable manifold follows from the invariance
of the diagonal and the fact that $O$ is a $1$-saddle if and only if
$\gamma>\gamma_1$.  For $\gamma<\gamma_1$, we will show independently that
there exist stationary points outside the diagonal. Note however that
relation~\eqref{ppA3} shows that these points must lie in a small
neighbourhood of the diagonal for $\gamma$ sufficiently close to
$\gamma_1$.
\end{proof}
 
Let us also point out that the growth of the potential away from the
diagonal can be controlled in the following way. 

\begin{prop}
\label{prop_ps2}
For any $x_\parallel\in\cD$ and $x_\perp$ orthogonal to the diagonal
$\cD$, the potential satisfies
\begin{equation}
\label{ps2}
V_\gamma(x_\parallel + x_\perp) \geqs V_\gamma(x_\parallel) + \frac12
\Bigpar{\frac\gamma{\gamma_1}-1}
\norm{x_\perp}^2\;.
\end{equation}
\end{prop}
\begin{proof}
Using~\eqref{ppA1} and the fact that $\Sigma x_\parallel=0$, we obtain 
for any $\lambda\in\R$ 
\begin{equation}
\label{ps2:1}
V_\gamma(x_\parallel+\lambda x_\perp) 
= -\frac12(\norm{x_\parallel}^2+\lambda^2\norm{x_\perp}^2)
+ \frac14 \sum_{i=1}^N (\brak{x_\parallel+\lambda x_\perp}_i)^4
+ \frac\gamma2 \lambda^2 \pscal{x_\perp}{\Sigma x_\perp}\;.
\end{equation}
The scalar product $\pscal{x_\perp}{\Sigma x_\perp}$ can be bounded below
by $\norm{x_\perp}^2/\gamma_1$. Moreover, applying Taylor's formula to
second order in $\lambda$, and using the fact that all components of
$x_\parallel$ are equal while the sum of the components of
$x_\perp$ vanishes,~\eqref{ps2} follows. 
\end{proof}

%%%%%%%%%%%%%%%%%%%%%%%%%%%%%%%%%%%%%%%%%%%%%%%%%%%%%%%%%%%%%%%%%%%%%%%%%%%

\section{Fourier Representation}
\label{sec_pf}

%%%%%%%%%%%%%%%%%%%%%%%%%%%%%%%%%%%%%%%%%%%%%%%%%%%%%%%%%%%%%%%%%%%%%%%%%%%

Let $\omega=\e^{2\pi\icx/N}$. The Fourier variables are defined by the
linear transformation
\begin{equation}
\label{pf1}
y_k = \frac1N \sum_{j\in\lattice} \cc{\omega}^{jk} x_j\;, 
\qquad\qquad
k \in \lattice^* = \Z/N\Z\;.
\end{equation}
The inverse transformation is given by 
\begin{equation}
\label{pf2}
x_j = \sum_{k \in \lattice^*} \omega^{jk} y_k\;,
\end{equation}
as a consequence of the fact that $\sum_{j=1}^N
\omega^{j(k-\ell)}=N\delta_{k\ell}$. Note that $y_k=\cc{y_{-k}}$, so that
we might use the real and imaginary parts of $y_k$ as dynamical
variables, instead of $y_k$ and $y_{-k}$. The following result is obtained
by a direct computation. 

\begin{prop}
\label{prop_pf1}
In Fourier variables, the equation of motion $\dot x=-\nabla V_\gamma(x)$
takes the form 
\begin{equation}
\label{pf3}
\dot y_k = \lambda_k y_k 
- \sum_{\substack{k_1,k_2,k_3 \in \lattice^* \\ k_1+k_2+k_3 = k}}
y_{k_1}y_{k_2}y_{k_3}\;,
\end{equation}
where the $\lambda_k$ are those defined in~\eqref{ppA2}. 
Furthermore, the potential is given in terms of Fourier variables by 
\begin{equation}
\label{pf4B}
\widehat V_\gamma(y) = -\frac N2 \sum_{k\in\Lambda^*} \lambda_k
\abs{y_k}^2 
+ \frac N4 
\sum_{\substack{k_1,k_2,k_3,k_4 \in \lattice^* \\ k_1+k_2+k_3+k_4 = 0}}
y_{k_1}y_{k_2}y_{k_3}y_{k_4}\;. 
\end{equation}
\end{prop} 

The effect of the symmetries on the Fourier variables is fully determined
by the action of three generators of the symmetry group $G$, as shown in
Table~\ref{table_Fourier_sym}. 

\begin{table}
\begin{center}
\begin{tabular}{|l|l|l|}
\hline
\vrule height 12pt depth 6pt width 0pt
$R$ & $x_j \mapsto x_{j+1}$ & $y_k \mapsto \omega^k y_k$ \\
\vrule height 8pt depth 6pt width 0pt
$RS=SR^{-1}$ & $x_j \mapsto x_{N-j}$ & $y_k \mapsto \cc{y_k} = y_{-k}$ \\
\vrule height 6pt depth 8pt width 0pt
$C$ & $x_j \mapsto -x_j$ & $y_k \mapsto -y_k$ \\
\hline 
\end{tabular} 
\end{center}
\caption[]
{Effect of a set of generators of the symmetry group $G_N$ on Fourier
variables.}
\label{table_Fourier_sym}
\end{table}

A particular advantage of the Fourier representation is that certain
invariant sets of phase space take a simple form in these variables. For
instance, for any $\ell$, 
\begin{align}
\label{pf55}
x_{\ell-j} &= x_j \;\forall j
&&\Rightarrow&
\cc{\omega}^{\ell k} \cc{y_k} &= y_k \;\forall k
&&\Rightarrow&
y_k &= \omega^{\ell k/2} r_k \;\forall k
\intertext{where the $r_k$ are all real. Similarly, for any $\ell$,}
x_{\ell-j} &= -x_j \;\forall j
&&\Rightarrow&
\cc{\omega}^{\ell k} \cc{y_k} &= -y_k \;\forall k
&&\Rightarrow&
y_k &= \icx \omega^{\ell k/2} r_k \;\forall k
\label{pf55v}
\end{align}
where the $r_k$ are all real. 

%%%%%%%%%%%%%%%%%%%%%%%%%%%%%%%%%%%%%%%%%%%%%%%%%%%%%%%%%%%%%%%%%%%%%%%%%%%

\subsection{The Case $N=2$}
\label{sec_f23}

For $N=2$, we have $\omega=-1$, and the Fourier variables are simply 
$y_0=(x_1+x_2)/2$, $y_1=(x_2-x_1)/2$. The equations~\eqref{pf3} become 
\begin{equation}
\label{f2_1}
\begin{split}
\dot y_0 &= y_0 \bigbrak{1-(y_0^2+3y_1^2)} \;, \\
\dot y_1 &= y_1 \bigbrak{\lambda_1-(3 y_0^2+y_1^2)} \;,
\end{split}
\end{equation}
with $\lambda_1=1-2\gamma$. The potential is given by 
\begin{equation}
\label{f2_1B}
\widehat V_\gamma(y) = -y_0^2 - \lambda_1 y_1^2 +
\frac12(y_0^4+6y_0^2y_1^2+y_1^4)\;.
\end{equation}

\begin{prop}
\label{prop_f2}
The bifurcation diagram for $N=2$ is the one given in \figref{fig_bif2},
and the value of the potential on the bifurcating branches is given by 
\begin{equation}
\label{f2_1C}
V_\gamma(A) = -\frac12\lambda_1^2\;,
\qquad
V_\gamma(Aa) = \frac1{16}(\lambda_1^2-6\lambda_1+1)\;.
\end{equation}
\end{prop}
\begin{proof}
In addition to the origin, there can be three types of stationary points:
\begin{itemiz}
\item	If $y_1=0$, $y_0\neq0$, then necessarily $y_0=\pm1$, yielding the
stationary points
$I^\pm$ in original variables.
\item	If $y_0=0$, $y_1\neq0$, there are two additional points $A, RA$,
given by
$y_1=\pm\sqrt{\lambda_1}$ whenever $\lambda_1>0$, i.e., $\gamma<1/2$. In
original variables, these have the expression
$(\pm\sqrt{\lambda_1},\mp\sqrt{\lambda_1})$, so that they 
have the $(x,-x)$-symmetry.
\item	If $y_0, y_1 \neq 0$, there are four additional points, given by
$8y_0^2=3\lambda_1-1$, $8y_1^2=3-\lambda_1$, provided $\lambda_1>1/3$,
i.e.,
$\gamma<1/3$. 
\end{itemiz}
It is straightforward to check the stability of these stationary points
from
the Jacobian matrix of~\eqref{f2_1}, and to compute the value of the
potential, using~\eqref{f2_1B}.
\end{proof}

%%%%%%%%%%%%%%%%%%%%%%%%%%%%%%%%%%%%%%%%%%%%%%%%%%%%%%%%%%%%%%%%%%%%%%%%%%%

\subsection{The Case $N=3$}
\label{sec_f3}

For $N=3$, we choose $\Lambda^*=\set{-1,0,1}$. The equations in Fourier
variables read 
\begin{equation}
\label{f3_1}
\begin{split}
\dot y_0 &= y_0 - (y_0^3 + y_1^3 + \yonebar^3 + 6y_0\abs{y_1}^2)\;, \\
\dot y_1 &= \lambda_1 y_1 - 3 (y_1\abs{y_1}^2 + y_0 \yonebar^2 + y_0^2
y_1)\;,
\end{split}
\end{equation}
with $\lambda_1 = 1 - \frac32\gamma$. 

\begin{prop}
\label{prop_f3}
The bifurcation diagram for $N=3$ is the one given in \figref{fig_bif3},
where the saddle--node bifurcations occur for 
\begin{equation}
\label{f3_1b}
\gamma = \gamma^\star(3) = \frac{\sqrt{3+2\sqrt3}-\sqrt3}3 \simeq
0.2701\dots\;.
\end{equation}
On the $1$-saddles, the potential has value
\begin{equation}
\label{f3_1c}
V_\gamma(A) = -\frac12\lambda_1^2\;.
\end{equation} 
\end{prop}
\begin{proof}
Using polar coordinates $y_1=r_1\e^{\icx\varphi_1}$, the
equations~\eqref{f3_1} become 
\begin{equation}
\label{f3_2}
\begin{split}
\dot y_0 &= y_0(1-y_0^2-6r_1^2) - 2r_1^3 \cos 3\varphi_1\;,\\
\dot r_1 &= r_1 \bigbrak{\lambda_1 - 3(r_1^2 + y_0 r_1 \cos 3\varphi_1 +
y_0^2)}\;,\\
\dot\varphi_1 &= 3y_0r_1\sin3\varphi_1\;.
\end{split}
\end{equation}
In addition to the origin, there can be three types of stationary points:
\begin{itemiz}
\item	If $r_1=0$, $y_0\neq0$, then necessarily $y_0=\pm1$, yielding the
points
$I^\pm$.
\item	If $y_0=0$, $r_1\neq0$, we obtain six stationary points given by 
$r_1=\sqrt{\lambda_1/3}$ and $\cos3\varphi_1=0$, provided $\lambda_1>0$,
that is, $\gamma<2/3$. These points have one of the symmetries $(x,-x,0)$,
$(x,0,-x)$ or $(0,x,-x)$. 
\item	If $\sin 3\varphi_1=0$, it is sufficient by symmetry to consider
the case $\varphi_1=0$ (i.e., $y_1$ real). These points have the
$(x,x,y)$-symmetry. Setting $y_0=u+v$, $r_1=u-v$, we find that stationary
points should satisfy the relations 
\begin{equation}
\label{f3_3}
\begin{split}
\lambda_1/3 &= 3u^2 + v^2 \;,\\
0 &= 24 u^2v + (1-\lambda_1)u + (1-\tfrac53\lambda_1)v\;.
\end{split}
\end{equation}
Taking the square of the second equation and eliminating $u$ yields a
cubic equation for $v^2$. In fact, the variable $z=v^2-(1+\lambda_1)/12$
satisfies 
\begin{equation}
\label{f3_4}
z^3 - \lambda z + \mu = 0\;, 
\qquad\qquad
\lambda = \frac{\lambda_1}{48}\;,
\qquad
\mu = \frac{1-3\lambda_1+6\lambda_1^2-2\lambda_1^3}{1728}\;.
\end{equation} 
This equation has three roots for $\lambda_1$ slightly smaller than $1$,
and
one root for $\lambda_1=0$. Bifurcations occur whenever the condition
$27\mu^2=4\lambda^3$ is fulfilled, which turns out to be equivalent to
$(1-\lambda_1)^2g(\lambda_1)=0$, where 
\begin{equation}
\label{f3_5}
g(\lambda_1) = 4\lambda_1^4 - 16\lambda_1^3 + 12\lambda_1^2 - 4\lambda_1 +
1\;.
\end{equation}
Since $g(0)=1$ and $g(1)=-3$, and it is easy to check that $g'<0$ on
$[0,1]$, there can be only one bifurcation point in this interval, whose
explicit value leads to the bifurcation value~\eqref{f3_1b}. 
\qed
\end{itemiz}
\renewcommand\qed{}
\end{proof}

%%%%%%%%%%%%%%%%%%%%%%%%%%%%%%%%%%%%%%%%%%%%%%%%%%%%%%%%%%%%%%%%%%%%%%%%%%%

\subsection{Centre-Manifold Analysis of the Desynchronisation Bifurcation}
\label{sec_cm}

Assume $N\geqs3$. 
We consider now the behaviour for $\gamma$ close to $\gamma_1$, i.e., for
$\lambda_1$ close to $0$. Setting $z=(y_0, y_2, y_{-2}, \dots)$, the
equation~\eqref{pf3} in Fourier variables is of the form 
\begin{equation}
\label{cm1}
\dot y_k = \lambda_k y_k + g_k(y_1, \yonebar, z)\;,
\end{equation}
where 
\begin{equation}
\label{cm2}
g_k(y_1, \yonebar, z) = - \sum_{\substack{k_1,k_2,k_3 \in \lattice^* \\
k_1+k_2+k_3 = k}}
y_{k_1}y_{k_2}y_{k_3}\;.
\end{equation}
For small $\lambda_1$, the system admits an invariant centre manifold of
equation 
\begin{equation}
\label{cm3}
y_k = h_k(y_1, \yonebar, \lambda_1)\;,\qquad k\neq \pm1\;,
\end{equation}
where the $h_k$ satisfy the partial differential equations 
\begin{equation}
\label{cm3B}
\lambda_k h_k + g_k(y_1, \yonebar, \set{h_j}_j) 
= \dpar{h_k}{y_1} \bigbrak{\lambda_1 y_1 + g_1(y_1, \yonebar,
\set{h_j}_j)} 
+ \dpar{h_k}{\yonebar} \cc{\bigbrak{\lambda_1 y_1 + g_1(y_1, \yonebar,
\set{h_j}_j)}}\;.
\end{equation}
We fix a cut-off order $K$. For our purposes, $K=2N$ will be sufficient. We
are looking for an expansion of the form 
\begin{equation}
\label{cm4}
h_k(y_1, \yonebar, \lambda_1) = \sum_{\substack{n,m\geqs0 \\ 3 \leqs n+m <
K}}
h^k_{nm}(\lambda_1) {y_1}^n \yonebar^m + \Order{\abs{y_1}^{K}}\;.
\end{equation}
First it is useful to examine the effect of symmetries on the
coefficients. 

\begin{lemma}
\label{lem_cm1}
The coefficients in the expansion of the centre manifold satisfy
\begin{itemiz}
\item	$h^k_{nm}(\lambda_1) \in \R$;
\item	$h^k_{nm}(\lambda_1) = h^{-k}_{mn}(\lambda_1)$;
\item	$h^k_{nm}(\lambda_1) = 0$ if $n-m \neq k \pmod N$;
\item	$h^k_{nm}(\lambda_1) = 0$ if $n+m$ is even.
\end{itemiz}
\end{lemma}
\begin{proof}
The centre manifold has the same symmetries as the equations~\eqref{pf3}. 
Thus  
\begin{itemiz}
\item	The $R$-symmetry requires
$\omega^k h_k(y_1, \yonebar, \lambda_1) 
= h_k(\omega y_1, \cc\omega\yonebar, \lambda_1)$, yielding 
the condition $(\omega^k-\omega^{n-m})h^k_{nm}(\lambda_1)=0$ so that
$h^k_{nm}(\lambda_1)$ vanishes unless $n-m = k \pmod N$;
\item	The $RS$-symmetry requires 
$\cc{h_k(y_1, \yonebar, \lambda_1)} = h_k(\yonebar, y_1, \lambda_1)$, and
yields the reality of the coefficients;
\item	The $C$-symmetry requires
$- h_k(y_1, \yonebar, \lambda_1) 
= h_k(- y_1, -\yonebar, \lambda_1)$, yielding 
the condition $((-1)^{n+m}+1)h^k_{nm}(\lambda_1)=0$, so that
$h^k_{nm}(\lambda_1)$ vanishes unless $n+m$ is odd;
\item	The condition 
$h_{-k}(y_1, \yonebar, \lambda_1) = \cc{h_k(y_1, \yonebar, \lambda_1)}$
yields the symmetry under permutation of $n$ and $m$.
\qed
\end{itemiz}
\renewcommand\qed{}
\end{proof}

From now on, we will write  $n-m\equiv k$ instead of $n-m=k\pmod N$. 
Lemma~\ref{lem_cm1} allows us to simplify the notation, setting 
\begin{equation}
\label{cm5}
h_k(y_1, \yonebar, \lambda_1) = \sum_{\substack{n,m\geqs0, \;3\leqs n+m < K
 \\ 
n-m\equiv k}}
h_{nm}(\lambda_1) {y_1}^n \yonebar^m + \Order{\abs{y_1}^{K}}\;.
\end{equation}
It is convenient to set $h_1(y_1, \yonebar,\lambda_1) = y_1$ and  
$h_{-1}(y_1, \yonebar,\lambda_1) = \yonebar$. Then~\eqref{cm5} holds for
$k=\pm1$ as well if we set $h_{nm}=\delta_{n1}\delta_{m0}$ whenever
$n-m\equiv1$, and $h_{nm}=\delta_{n0}\delta_{m1}$ whenever $n-m\equiv-1$. 
The equation on the centre manifold can be written as 
\begin{align}
\nonumber
\dot y_1 &= \lambda_1 y_1 - \sum_{k_1+k_2+k_3 = 1}
h_{k_1}(y_1,\yonebar,\lambda_1) 
h_{k_2}(y_1,\yonebar,\lambda_1)
h_{k_3}(y_1,\yonebar,\lambda_1) \\
&= 
\lambda_1 y_1 - \sum_{\substack{n,m\geqs0, \;3\leqs n+m < K  \\ n-m\equiv
1}} c_{nm}(\lambda_1) y_1^n \yonebar^m + \Order{\abs{y_1}^K}\;,
\label{cm6}
\end{align}
where 
\begin{equation}
\label{cm7}
c_{nm}(\lambda_1) = 
\sum_{\substack{n_i\geqs 0 \colon n_1+n_2+n_3=n \\ 
m_i\geqs 0 \colon m_1+m_2+m_3=m}} 
h_{n_1m_1}(\lambda_1) h_{n_2m_2}(\lambda_1) h_{n_3m_3}(\lambda_1)
\in\R\;.
\end{equation}
Note that $c_{nm}(\lambda_1) = 0$ whenever $n+m$ is even. 
In polar coordinates $y_1 = r_1 \e^{\icx\varphi_1}$, Equation~\eqref{cm6}
becomes 
\begin{equation}
\label{cm8}
\begin{split}
\dot r_1 &= \lambda_1 r_1 - 
\sum_{\substack{n,m\geqs0, \;3\leqs n+m < K  \\ n-m\equiv 1}} 
c_{nm}(\lambda_1) r_1^{n+m} \cos\bigpar{(n-m-1)\varphi_1} +
\Order{r_1^K}\;, \\
\dot\varphi_1 &= 
-\sum_{\substack{n,m\geqs0, \;3\leqs n+m < K  \\ n-m\equiv 1}} 
c_{nm}(\lambda_1) r_1^{n+m-1} \sin\bigpar{(n-m-1)\varphi_1} +
\Order{r_1^{K-1}}\;.
\end{split}
\end{equation}
In general, $c_{21} = 3h_{10}^2h_{01} = 3$ is the only term contributing to
the third-order term of $\dot r_1$. The only exception is
the case $N=4$, in which $c_{03} = h_{10}^3 = 1$ also contributes to the
lowest order, yielding 
\begin{equation}
\label{cm9}
\dot r_1 = 
\begin{cases}
\lambda_1 r_1 - 3 r_1^3 + \Order{r_1^5}
& \text{if $N\neq 4$\;,} \\
\lambda_1 r_1 - (3+\cos4\varphi_1)r_1^3 + \Order{r_1^5}
& \text{if $N=4$\;.}
\end{cases}
\end{equation}
This shows that all stationary points bifurcating from the origin  lie at a
distance of order $\sqrt{\lambda_1}$ from it: They satisfy  $r_1 =
\sqrt{\lambda_1/(3+\cos4\varphi_1)} + \Order{\lambda_1^{3/2}}$ in the case
$N=4$, and $r_1 = \sqrt{\lambda_1/3} + \Order{\lambda_1^{3/2}}$ otherwise. 

The terms with $n-m=1$ do not contribute to the angular derivative
$\dot\varphi_1$. In the particular case $N=4$, we have 
\begin{equation}
\label{cm10}
\dot\varphi_1 = \sin(4\varphi_1) r_1^2 + \Order{r_1^4}\;,
\end{equation}
yielding $8$ stationary points, of alternating stability. Otherwise, we
have
to distinguish between two cases:
\begin{itemiz}
\item	If $N$ is even, the lowest-order coefficient contributing to
$\dot\varphi_1$ is $c_{0,N-1}$, giving 
\begin{equation}
\label{cm11}
\dot\varphi_1 = c_{0,N-1}(\lambda_1) r_1^{N-2} \sin(N\varphi_1) +
\Order{r_1^{N-1}}\;.
\end{equation}
Thus if we prove that $c_{0,N-1}(\lambda_1)\neq 0$, we will have obtained
the existence
of exactly $2N$ stationary points, of alternating stability, 
bifurcating from the origin.

\item	If $N$ is odd, then $n+m$ is even whenever $n-m = \pm N+1$, which
implies by Lemma~\ref{lem_cm1} that $c_{nm}=0$ for these $(n,m)$. The
lowest-order coefficient contributing to $\dot\varphi_1$ is thus
$c_{0,2N-1}$, giving
\begin{equation}
\label{cm12}
\dot\varphi_1 = c_{0,2N-1}(\lambda_1) r_1^{2N-2} \sin(2N\varphi_1) +
\Order{r_1^{2N-1}}\;.
\end{equation}
Thus if we prove that $c_{0,2N-1}(\lambda_1)\neq 0$, we will have obtained
the existence
of exactly $4N$ stationary points, of alternating stability, 
bifurcating from the origin.
\end{itemiz}

\begin{remark}
\label{rem_cm0}
Let $r_1$ and $r'_1$ be solutions of the equation $\dot r_1=0$ obtained for
two different values of $\varphi_1$, say $\varphi_1=0$ and
$\varphi_1=\pi/N$. For even $N$, using the fact that the terms up to order
$N-2$ in the first equation in~\eqref{cm8} do not depend on $\varphi_1$,
one can see that $r'_1-r_1=\Order{\lambda_1^{(N-3)/2}}$. For odd $N$, one 
obtains in a similar way $r'_1-r_1=\Order{\lambda_1^{(2N-3)/2}}$. 
\end{remark}

In order to compute sign of the coefficients $c_{0,N-1}$ or $c_{0,2N-1}$,
we need at
least to know the coefficients $h_{0m}$ for odd $m$ up to $N-3$ or $2N-3$,
respectively. By continuity, however, it is sufficient to compute them for
$\lambda_1=0$. We henceforth set $h_{nm}=h_{nm}(0)$.

\begin{lemma}
\label{lem_cm0}
For all odd $m>0$ such that $m\not\equiv\pm1$, $h_{0m}=h_{m0}$ satisfies
\begin{align}
\nonumber
\lambda_m h_{0m}
={}& 
\sum_{m_i\geqs 0 \colon m_1+m_2+m_3=m}
h_{0m_1} h_{0m_2} h_{0m_3} \\
\nonumber
&{}- \sum_{\substack{v \geqs 0 \colon v \equiv m+1 \\
m_i\geqs 0 \colon m_1+m_2+m_3+v=m}}
h_{0m_1} h_{0m_2} h_{0m_3} h_{1v} \\
&{}- \sum_{\substack{v > 0 \colon v \equiv m \\
n_i\geqs 0 \colon n_1+n_2+n_3+v=m+1}}
v h_{n_10} h_{n_20} h_{n_30} h_{0v}\;.
\label{cm15}
\end{align}
Furthermore, if either $N$ is even and $1\leqs m\leqs N-3$, or $N$ is odd
and $1\leqs m\leqs 2N-3$, then 
\begin{equation}
\label{cm17}
\lambda_m h_{0m}
= \sum_{m_i\geqs 0 \colon m_1+m_2+m_3=m}
h_{0m_1} h_{0m_2} h_{0m_3}\;.
\end{equation}
\end{lemma}
\begin{proof}
By invariance of the centre manifold, $h_k(y_1, \yonebar,0)$ has to satisfy
the equation 
\begin{equation}
\label{cm13}
\lambda_k h_k = -g_k(y_1, \yonebar, \set{h_j}_j) 
+ \dpar{h_k}{y_1} g_1(y_1, \yonebar, \set{h_j}_j) 
+ \dpar{h_k}{\yonebar} \cc{g_1(y_1, \yonebar, \set{h_j}_j)}\;.
\end{equation}
Plugging in the series~\eqref{cm2} of $g_k$ and~\eqref{cm5} of $h_k$, this
can be seen to be equivalent to 
\begin{align}
\nonumber
\lambda_{n-m} h_{nm}
={}& 
\sum_{\substack{n_i\geqs 0 \colon n_1+n_2+n_3=n \\ m_i\geqs 0 \colon
m_1+m_2+m_3=m}}
h_{n_1m_1} h_{n_2m_2} h_{n_3m_3} \\
\nonumber
&{}- \sum_{\substack{u, v \geqs 0 \colon u - v \equiv n - m \\
n_i\geqs 0 \colon n_1+n_2+n_3+u=n+1 \\ m_i\geqs 0 \colon m_1+m_2+m_3+v=m}}
u h_{n_1m_1} h_{n_2m_2} h_{n_3m_3} h_{uv} \\
&{}- \sum_{\substack{u, v \geqs 0 \colon u - v \equiv n - m \\
n_i\geqs 0 \colon n_1+n_2+n_3+v=m+1 \\ m_i\geqs 0 \colon m_1+m_2+m_3+u=n}}
v h_{n_1m_1} h_{n_2m_2} h_{n_3m_3} h_{uv}\;.
\label{cm14}
\end{align}
In the special case $n=0$, the second sum vanishes unless $u=1$, and
$n_1=n_2=n_3=0$. In the third sum, we must have $v>0$ and
$u=m_1=m_2=m_3=0$.
Finally, the fact that $\lambda_{-m}=\lambda_m$ yields~\eqref{cm15}. 

Assume now that $N$ is even and $m\leqs N-3$. In the second sum, $v$
cannot exceed $m$, but then the condition $v\equiv m+1$ would require
$m\geqs N-1$. Thus the second sum vanishes. In the third sum, $v$ cannot
exceed $m+1$, and thus $v=m$. However, in that case $n_1+n_2+n_3=1$, so
that two $n_i$ must be zero. Since $h_{00}=0$, the third sum vanishes as
well. 

If $N$ is odd and $m\leqs 2N-3$, then the second sum in~\eqref{cm15} allows
for $v=m+1-N$. Then, however, we would have $m_1+m_2+m_3=N-1$, which is
even. Thus at least one of the $m_i$ is even, yielding a vanishing summand.
The third sum in~\eqref{cm15} allows for $v=m$ and $v=m-N$. In the first
case, however, the summand vanishes for the same reason as before, while in
the second case, we would have $n_1+n_2+n_3=N+1$, which is even. Thus at
least one of the $n_i$ is even, yielding again a vanishing summand.  
\end{proof}

\begin{prop}
\label{prop_cm1}
If $N$ is even, then 
\begin{equation}
\label{cm16}
\begin{cases}
c_{0,N-1} > 0 & \text{if $N\in 4\N$\;,} \\
c_{0,N-1} < 0 & \text{if $N\in 4\N + 2$\;.}
\end{cases}
\end{equation}
As a consequence, for $\lambda_1 > 0$ sufficiently small, the system admits
exactly $2N$ stationary points on the centre manifold. The points with
$\varphi_1 = 2 k\pi/N$ have one stable and one unstable direction if
$N\in 4\N$, and two stable directions if $N\in 4\N + 2$, and vice versa for
the points with $\varphi_1 = (2k+1)\pi/N$. 
\end{prop}
\begin{proof}
By~\eqref{cm7}, it is sufficient to compute $h_{0m}=h_{m0}$ for odd $m$
between $1$ and $N-3$.  Recall that $\lambda_0=1$, and
$\lambda_k=\lambda_{-k}<0$ for $3\leqs k\leqs N-2$. Using $h_{01}=1$ as
starting point, it is easy to show by induction that $\sign(h_{0,2\ell+1})
= (-1)^\ell$, because all summands have the same sign at each iteration.
Likewise, all summands of $c_{0,N-1}$ have sign $(-1)^{N/2}$. 
\end{proof}

\begin{proof}[{\sc Proof of Theorem~\ref{thm_desync1}}]
We consider the case $N=4L$, the proof being similar for $N=4L+2$.
\begin{itemiz}
\item	First note that the set $\cB=\setsuch{y}{y_k\in\R\;\forall k}$ is
invariant under the dynamics (it corresponds to $x_{N-j}=x_j\;\forall j$). 
The intersection of $\cB$ with the centre manifold is one-dimensional, and
can be parametrised by $r_1\in\R$ (while $\varphi_1=0$). Since $\dot r_1 =
\lambda_1 r_1 - 3r_1^3 + \Order{r_1^4}$, the system admits at least three
stationary points $O$ and $\pm B'$ in $\cB$ for small positive $\lambda_1$.
The stationary point $B'$ is stable in the $r_1$-direction, and unstable in
the $\varphi_1$-direction. Since there are $N-3$ stable and $1$ unstable
directions transversal to the centre manifold, $B'$ is a $2$-saddle. The
same holds for the cyclic permutations $RB',\dots,R^{N-1}B'$, which
correspond to $\varphi_k=2k\pi/N$, and thus lie in the fixed-point sets of
conjugate
symmetry groups. Applying the inverse Fourier transformation~\eqref{pf2},
we find that the
coordinates of $B'$ in $\cX$ satisfy 
\begin{equation}
\label{cm17:1}
B'_j = \omega^j y_1 + \omega^{-j} \yonebar + \Order{\lambda_1}
= 2\cos \biggpar{\frac{2\pi j}N} \sqrt{\frac{\lambda_1}3} +
\Order{\lambda_1}\;.
\end{equation}
Setting $B=R^{-N/4}B'$ yields the
expression~\eqref{desybif3} for the coordinates. 

\item	A similar argument shows the existence of $N$ stationary points
$A,RA,\dots,R^{N-1}A$, corresponding to $\varphi_1=(2k+1)\pi/N$, which are
$1$-saddles because they are stable in the $\varphi_1$-direction. In $\cX$,
their coordinates satisfy one of the symmetries $x_{n_0-j}=-x_j$,
$n_0=1,\dots,N$. 

\item	For $\lambda_1$ small enough, the equation $\dot\varphi_1=0$ admits
exactly $2N$ solutions, so that there are no further stationary points on
the centre manifold. Proposition~\ref{prop_ps1} shows that for any
$\eta>0$,
we can find a $\delta>0$ such that if $\gamma>\gamma_1-\delta$, there can
be
no stationary points outside an $\eta$-neighbourhood of the diagonal.
Together with a local analysis near the diagonal, and the fact that the
centre manifold is locally repulsive, this proves that there are exactly
$2N+3$ stationary points. 

\item	Consider now the set 
\begin{equation}
\label{cm17:2}
\cA_+ = \setsuch{x\in\cX}{x_j=-x_{N+1-j}=x_{N/2+1-j} \;\forall j,\;
x_1,\dots,x_L>0}\;.
\end{equation}
We claim that this set is positively invariant under the flow of $\dot
x=-\nabla V_\gamma(x)$. Without the condition $x_1,\dots,x_L>0$, the
invariance
follows from equivariance. Now if $x_j=0$ for some $j$ while the other
$x_i$
are positive, one easily sees that $-\sdpar
{V_\gamma}{x_j}=(x_{j-1}+x_{j+1})\gamma/2$ is positive, showing the
invariance of
$\cA_+$. 
Since the potential is increasing at infinity, there must be at
least one stationary point in $\cA_+$, which we denote $A(\gamma)$. As
$\gamma\to 0$, the only stationary point in $\cA_+$ is the point
$(1,\dots,1,-1,\dots,-1)$. We proceed similarly for $B(\gamma)$. 

\item	Finally, the value of the potential at the stationary points can be
computed with the help of the expression~\eqref{pf4B} for the potential in
Fourier variables. The only terms contributing to leading order in
$\lambda_1$ are the term $\lambda_1\abs{y_1}^2$ in the first sum, and terms
of the form $y_1^2y_{-1}^2$ in the second sum. The remaining terms are of
smaller order. The relation on the difference $V_\gamma(B)-V_\gamma(A)$ is
a consequence of Remark~\ref{rem_cm0}.  This proves the theorem.  \qed 
\end{itemiz} \renewcommand\qed{} \end{proof}

In the case of odd $N$, the situation is more difficult, because not all
summands in the recursion defining the $h_{0m}$ are of the same sign. 
A partial result is 

\begin{lemma}
\label{lem_cm_odd}
If $N$ is odd, then  
\begin{equation}
\label{cm18}
\sign(h_{0,2\ell+1}) = (-1)^\ell
\qquad
\text{for $\ell = 0,1,\dots,\frac{N-1}2-1$\;.}
\end{equation}
\end{lemma}
\begin{proof}
As before, using the fact that
$\lambda_k=\lambda_{-k}<0$ for $1\leqs k\leqs N-1$. 
\end{proof}

\begin{conjecture}
\label{con_cm_odd}
For $\ell = (N-1)/2,\dots,N-2$, $h_{0m}$ has sign $(-1)^{\ell+1}$, and 
\begin{equation}
\label{cm20}
c_{0,2N-1} > 0\;.
\end{equation}
\end{conjecture}

Numerically, we have checked the validity of this conjecture for all odd
$N$
up to $101$. 

\begin{proof}[{\sc Proof of Theorem~\ref{thm_desybif2}}] 
The proof is similar to the above proof of Theorem~\ref{thm_desync1},
without using any information on the sign of $c_{0,2N-1}$. 
\end{proof}

Conjecture~\ref{con_desync2} then follows from
Conjecture~\ref{con_cm_odd} by including the information on the sign of
$c_{0,2N-1}$ in the proof.

%%%%%%%%%%%%%%%%%%%%%%%%%%%%%%%%%%%%%%%%%%%%%%%%%%%%%%%%%%%%%%%%%%%%%%%%%%%

\newpage
\appendix
\section{Small Coupling and Symbolic Dynamics}
\label{app_horse}

In this appendix we sketch the proof of Proposition~\ref{prop_sl1} on the
continuation of equilibrium points from the uncoupled limit.  The equation
$f(x_n)+\frac\gamma2(x_{n+1}-2x_n+x_{n-1})=0$ satisfied by the stationary
points of the potential $V_\gamma$ can be rewritten as
$(x_{n+1},y_{n+1})=H(x_n,y_n)$ where $H$ is the map of the plane defined
by 
\begin{equation}
H(x,y)=\Bigpar{2x-\frac{2}{\gamma}f(x)-y,x}\;.
\end{equation}
This map is invertible and its inverse can be simply obtained as
$H^{-1}=S\circ H\circ S$ where $S(x,y)=(y,x)$. In addition, we have
$H(-x,-y)=-H(x,y)$.

We proceed to a similar construction as in~\cite{Keener87} to show that,
when $\gamma\leqs1/4$, the map $H$ has an invariant set on which it is
conjugated to the full shift with $3$ symbols. This implies that the
map admits at least $3^N$ points of (not necessarily minimal) period $N$
for any $N$. Since the potential $V_\gamma$ is a polynomial of degree $4$
in $N$ variables, it can have at most $3^N$ stationary points, and thus we
will have obtained the exact number of stationary points. 

\begin{figure}
\centerline{\includegraphics*[clip=true,width=130mm]{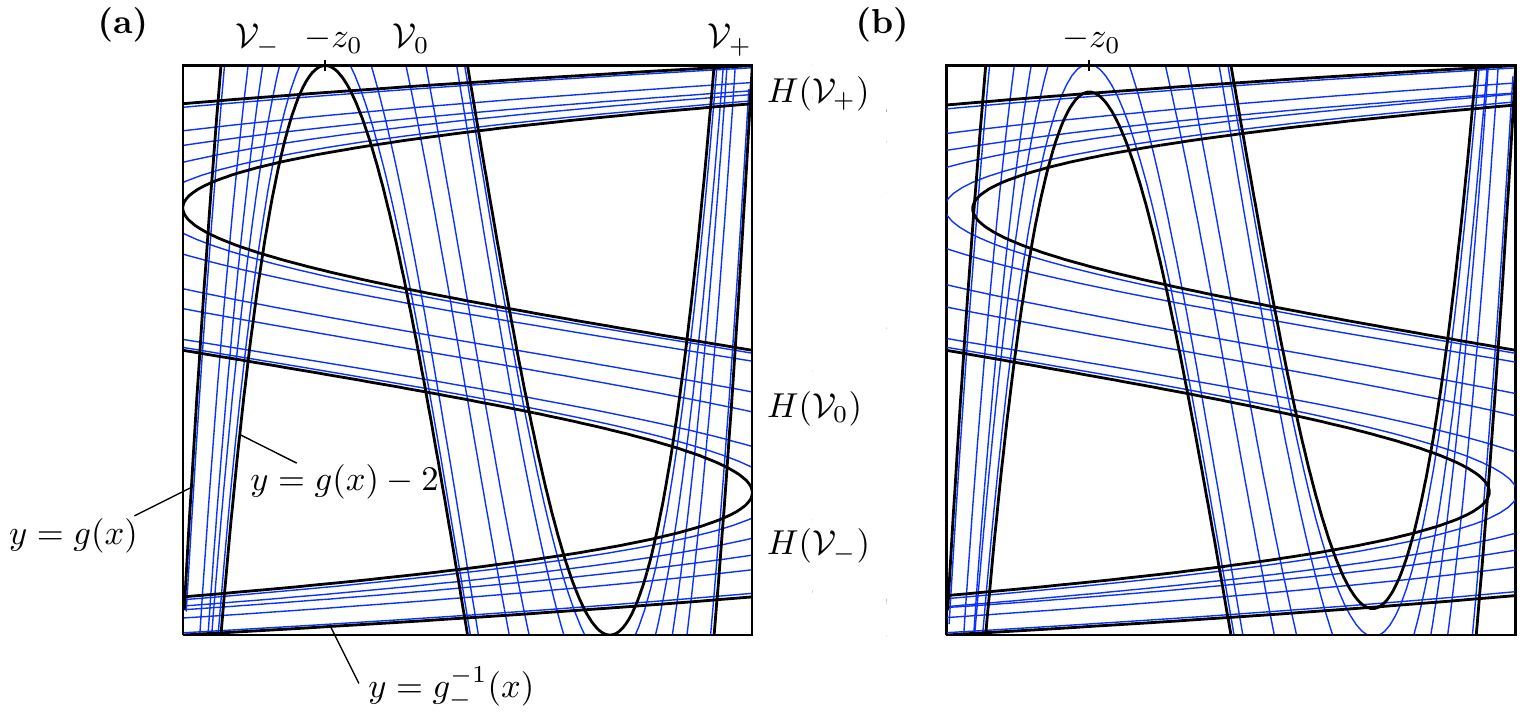}}
 \figtext{ }
 \caption[]
 {Images of the boundary of the square $[-1,1]\times[-1,1]$ under the map
 $H$, shown {\bf (a)} for $\gamma=1/4$, and {\bf (b)} for
$\gamma=0.258\dots$
 The \lq\lq vertical\rq\rq\ strips $\cV_-$, $\cV_0$ and $\cV_+$ are mapped
on
 \lq\lq horizontal\rq\rq\ strips $H(\cV_-)$, $H(\cV_0)$ and $H(\cV_+)$.
Their
 intersections define the elements of the first-level partition, for
 instance $I_{-.-}=\cV_-\cap H(\cV_-)$ is the lower-left squarish domain.
 Intersections of the light curves, obtained by iterating the map twice,
 define the elements $I_{\omega_{-1}\omega_0.\omega_1\omega_2}$ 
 of the second-level partition.}
\label{fig_horse}
\end{figure}

We start by defining a collection of strips in the square $[-1,1]^2$ with
good dynamical properties.
The map $x\mapsto g(x)\defby2x-2\gamma^{-1}f(x)+1$ leaves $-1$ invariant
and is strictly increasing on $[-1,-z_0]$ where
$z_0\defby\sqrt{(1-\gamma)/3}$. If in addition, the condition $g(-z_0)\geqs
1$ holds, then $g$ is a bijection on $[-1,g_-^{-1}(1)]$ where $g_-^{-1}$ is
the inverse of $g$ with range $[-1,-z_0]$ and we have 
\begin{equation}
H\Bigpar{\bigsetsuch{(x,-1)}{g_-^{-1}(-1)\leqs x\leqs g_-^{-1}(1)}}
=\bigsetsuch{(x,g_-^{-1}(x))}{-1\leqs x\leqs 1}\;.
\end{equation}
Similarly, provided that the condition $g(-z_0)\geqs 3$ is satisfied, we
have 
\begin{equation}
H\Bigpar{\bigsetsuch{(x,1)}{g_-^{-1}(1)\leqs x\leqs g_-^{-1}(3)}}=
\bigsetsuch{(x,g_-^{-1}(x+2))}{-1\leqs x\leqs 1}\;.
\end{equation}
The condition $g(-z_0)\geqs 3$ is indeed equivalent to
$\gamma\leqs1/4$. By monotonicity, it results that under the
condition $\gamma\leqs1/4$ the \lq\lq vertical\rq\rq\ strip 
\begin{equation}
\cV_- \defby 
\bigsetsuch{(x,y)}{g_-^{-1}(y)\leqs x\leqs g_-^{-1}(y+2),\ -1\leqs y\leqs
1}\;,
\end{equation}
is mapped onto the \lq\lq horizontal\rq\rq\ strip $S(\cV_-)$.  By symmetry
we
also have that $H(\cV_+)=S(\cV_+)$ where $\cV_+=-\cV_-$ is the opposite
vertical
strip.

Furthermore, $g$ is strictly decreasing on $[-z_0,z_0]$ and the conditions
$g(-z_0)\geqs 1$ and $g(z_0)\leqs -1$ imply that this map is a bijection on
$[g_+^{-1}(1),g_+^{-1}(-1)]$ where $g_+^{-1}$ is the inverse of $g$ with
range $[z_0,-z_0]$. The inverse $g_+^{-1}$ has the symmetry
$g_+^{-1}(2-x)=-g_+^{-1}(x)$ and the condition $g(z_0)\leqs -1$ is
equivalent to $g(-z_0)\geqs 3$. Similarly as before this implies that under
the condition $\gamma\leqs1/4$, the symmetric vertical strip 
\begin{equation}
\cV_0\defby 
\bigsetsuch{(x,y)}{g_+^{-1}(y+2)\leqs x\leqs g_+^{-1}(y),\ -1\leqs y\leqs
1}\;,
\end{equation}
is mapped onto the horizontal strip $S(\cV_0)$.

With these strips provided, the rest of the construction is standard. Given
any finite word $\omega_{-n}\cdots\omega_{n+1}\in\{-,0,+\}^{2(n+1)}$
($n\geqs 0$) consider the set
\begin{equation}
\label{horse6}
I_{\omega_{-n}\cdots\omega_0.\omega_1\cdots\omega_{n+1}} 
= \bigcap_{k=-n}^{n+1} H^k(\cV_{\omega_k})\;.
\end{equation}
These sets are all non-empty by the properties of the strips above. Then,
for any infinite sequence $\{\omega_n\}\in\{-,0,+\}^\Z$, the intersection
$\bigcap_{k\in\Z} H^k(\cV_{\omega_k})$ is also nonempty. Since the
sets $I_{\omega_{-n}\cdots\omega_{n+1}}$ are pairwise disjoint, we conclude
that the map $H$ has at least (and thus exactly) $3^\Z$ points for which
the forward and backward orbit is contained in $[-1,1]^2$. In particular it
has exactly $3^N$ points of period $N$ for any
$N$. This shows that $\gamma^\star(N)\geqs 1/4$ for any $N\geqs 2$.  

The condition $g(-z_0)\geqs 3$ which ensures that the strips have the
desired properties is not optimal and can be improved. To see this, we
consider the strips $\cV'_{\omega_k}$ instead of $\cV_{\omega_k}$ in
\eqref{horse6} where $\cV'_-$ is defined by 
\begin{equation}
\cV'_-\defby\bigsetsuch{(x,y)}{-1\leqs x\leqs -z_0,\ g(x)-2\leqs y\leqs
g(x)}\;,
\end{equation}
where $\cV'_+=-\cV'_-$ and where 
\begin{equation}
\cV'_0\defby\bigsetsuch{(x,y)}{-z_0\leqs x\leqs z_0,\ g(x)-2\leqs y\leqs
g(x)}\;,
\end{equation}
Note that we have $\cV'_{\omega_k}\supset \cV_{\omega_k}$. 

Similarly as before, a sufficient condition for the existence of $3^\Z$
points with bounded orbit under $H$ is that  the nine basic sets
$I_{\omega_0.\omega_1}$ are pairwise disjoint. By symmetries, the latter
is equivalent to the condition that the right boundary of $\cV'_-$ crosses
the upper horizontal strip $S(\cV'_+)$. A simple analysis shows that this
is
equivalent to the inequality (weaker than $g(-z_0)\geqs 3$)
\begin{equation}
g(-z_0)-2\geqs g_-^{-1}(2-z_0)
\end{equation}
which holds provided $\gamma\leqs 0.258\dots$, i.e., we have
$\gamma^\star(N)\geqs 0.258\dots$ for any $N\geqs 2$. 

As $\gamma$ increases, some symbolic sequences will be pruned, resulting in
saddle-node bifurcations of stationary points. The first configurations to
disappear are asymmetric ones, such as $(1,1,\dots,1,-1)$ and
$(1,1,\dots,1,0)$ and their images under the symmetry group~$G$. 

%%%%%%%%%%%%%%%%%%%%%%%%%%%%%%%%%%%%%%%%%%%%%%%%%%%%%%%%%%%%%%%%%%%%%%%%%%%

\newpage
\section{The case $N=4$}
\label{app_N4}

In this appendix, we briefly describe how we obtained the bifurcation
diagram of~\figref{fig_bif4} for the case $N=4$. 

We start by determining the isotropy groups of the stationary points
present in the uncoupled case $\gamma=0$, and their fixed-point sets. 
They are of $10$ different types, $8$ of which we show in
Table~\ref{table_N=4}. The two types we do not show are the fixed-point set
$\set{0,0,0,0}$, which obviously contains only the origin, and the orbits
labelled $Aa\alpha$, which have no symmetry, and thus a fixed-point set
equal
to $\cX$. 

\begin{table}[t]
\begin{center}
\begin{tabular}{|l|l|l|}
\hline
\vrule height 12pt depth 6pt width 0pt
Orbit label & Fixed-point set & Reduced dynamics \\
\hline
\vrule height 12pt depth 6pt width 0pt
$I^\pm$ & $(x,x,x,x)$ & $\dot x = x-x^3$ \\
\vrule height 6pt depth 6pt width 0pt
$A$ & $(x,x,-x,-x)$ & $\dot x = (1-\gamma)x-x^3$ \\
\vrule height 6pt depth 6pt width 0pt
$B$ & $(x,0,-x,0)$ & $\dot x = (1-\gamma)x-x^3$ \\
\vrule height 6pt depth 8pt width 0pt
$A^{(2)}$ & $(x,-x,x,-x)$ & $\dot x = (1-2\gamma)x-x^3$ \\
\hline 
\vrule height 12pt depth 6pt width 0pt
$Aa$ & $(x,x,y,y)$ 	& $\dot x = (1-\frac12\gamma)x + \frac12\gamma y -
x^3$ \\
\vrule height 6pt depth 8pt width 0pt
 & 		& $\dot y = \frac12\gamma x + (1-\frac12\gamma)y - y^3$ \\
\vrule height 8pt depth 6pt width 0pt
 & $(x,y,x,y)$  & $\dot x = (1-\gamma)x + \gamma y - x^3$ \\
\vrule height 6pt depth 8pt width 0pt
 & 		  & $\dot y = \gamma x + (1-\gamma)y - y^3$ \\
\vrule height 8pt depth 6pt width 0pt
 & $(x,-x,y,-y)$ & $\dot x = (1-\frac32\gamma)x - \frac12\gamma y - x^3$ \\
\vrule height 6pt depth 8pt width 0pt
 & 		   & $\dot y = -\frac12\gamma x + (1-\frac32\gamma)y - y^3$
\\ 
\hline 
\vrule height 12pt depth 6pt width 0pt
$\partial a, \partial b, \dots$ & $(x,y,x,z)$ 	& $\dot x = (1-\gamma)x +
\frac12\gamma y + \frac12\gamma z- x^3$ \\
\vrule height 8pt depth 6pt width 0pt
 & 		& $\dot y = \gamma x + (1-\gamma)y - y^3$ \\
\vrule height 6pt depth 7pt width 0pt
 & 		& $\dot z = \gamma x + (1-\gamma)z - z^3$ \\ 
\hline 
\end{tabular} 
\end{center}
\caption[]
{Orbits, corresponding fixed-point sets and reduced dynamics in the case
$N=4$.}
\label{table_N=4}
\end{table}

Table~\ref{table_N=4} also shows the form taken by the equation $\dot
x=-\nabla V_\gamma(x)$ when restricted to the fixed-point set. These
equations are
then analysed in order to determine their equilibrium points. The analysis
of the dynamics in one-dimensional fixed-point sets is straightforward,
while the two-dimensional equations require a bit more work, but do not
present particular difficulties. Let us briefly describe the case of the
three-dimensional fixed-point set $(x,y,x,z)$.  In rotated variables
$u=(y+z)/2$, $v=(y-z)/2$, the equations become
\begin{equation}
\label{appB1}
\begin{split}
\dot x &= (1-\gamma) x + \gamma u - x^3\;,\\
\dot u &= \gamma x + (1-\gamma) u - u(u^2+3v^2)\;,\\
\dot v &= (1-\gamma) v - v(3u^2+v^2)\;.
\end{split}
\end{equation}
Equilibrium points with $v=0$ correspond to the fixed-point set
$(x,y,x,y)$,
and have already been analysed. We thus assume $v\neq0$, which implies
$3u^2+v^2=1-\gamma$ and allows to eliminate $v$. We are left with two
equations for the stationary points, namely   
\begin{equation}
\label{appB2}
\begin{split}
(1-\gamma) x + \gamma u - x^3 &= 0\;,\\
\gamma x - 2(1-\gamma) u + 8u^3 &= 0\;.
\end{split}
\end{equation}
If $x=0$, then necessarily $u=0$ and we are left with the already studied
points of the form $(0,y,0,-y)$. We thus assume $x\neq0$ and eliminate $u$
from the system. It is convenient to introduce the variable
$w=2(1-\gamma-x^2)$ instead of $x$. Then $u=-wx/2\gamma$, and $w$ has to
satisfy the fourth-order equation 
\begin{equation}
\label{appB3}
w^4 - 2(1-\gamma)w^3 + 2\gamma^2(1-\gamma)w + 2\gamma^4 = 0\;.
\end{equation}
The condition $v^2\geqs 0$ yields the additional requirement 
$w^3 - 2(1-\gamma)w^2 + \frac83\gamma^2(1-\gamma) \geqs 0$. Together
with~\eqref{appB3}, this can be seen to be equivalent to the condition
\begin{equation}
\label{appB4}
\frac{3\gamma^2}w \leqs 1 - \gamma\;.
\end{equation}
Numerically, one observes the following properties:
\begin{itemiz}
\item	If $0 < \gamma < 0.2684\dots$, Equation~\eqref{appB3} has four
distinct
real roots, two of them being positive and two negative. 
\item	For $0.2684\dotsc < \gamma < 0.4004\dots$, Equation~\eqref{appB3}
has
two real roots, which are both positive. 
\end{itemiz}
The negative roots always fulfil Condition~\eqref{appB4}.
Examining~\eqref{appB3} perturbatively in $\gamma$, one finds that they are
of the form $w=-\gamma+\Order{\gamma^2}$ and
$w=-\gamma^2+\Order{\gamma^3}$, and correspond to the branches labelled
$\partial a$ and $\partial b$ in the bifurcation diagram
of~\figref{fig_bif4}. 

The positive roots, on the other hand, do not always fulfil
Condition~\eqref{appB4}. In fact, one solution exists for
$0<\gamma<(3\sqrt2-2)/7$, while the other one bifurcates with the
$A^{(2)}$-branch for $\gamma=2/5$. A perturbative analysis for small
$\gamma$ shows that they converge, respectively, to the points $(1,0,1,-1)$
and $(0,1,0,0)$ as $\gamma\to0$. 

Finally, one can also compute the determinant of the Hessian of the
potential around the stationary points of the form $(x,y,x,z)$. This
determinant can be put into the form 
\begin{equation}
\label{appB5}
\frac1{2w^2} \bigbrak{3w-4(1-\gamma)} \bigbrak{(1-\gamma)w - 3\gamma^2} 
\bigbrak{3(1-\gamma)w^2 + 4(1-2\gamma)w + 6\gamma^2(1-\gamma)}\;. 
\end{equation}
The roots of this expression correspond to bifurcation points, and can be
combined with~\eqref{appB3} to check that all bifurcation points have been
found (and yields more precise estimates of the bifurcation values). This
expression can also be used to check the index of the stationary points. 

This analysis is not completely rigorous because we have not analysed in
detail the branch labelled $Aa\alpha$, which has no symmetry at all. The
bifurcation diagram is based on the fact that the only bifurcation of
already known branches producing stationary points without symmetry is the
bifurcation of the $Aa$-branch at $\gamma=1/3$. A local analysis shows that
this is indeed a pitchfork bifurcation, producing the right number of new
stationary points.

%%%%%%%%%%%%%%%%%%%%%%%%%%%%%%%%%%%%%%%%%%%%%%%%%%%%%%%%%%%%%%%%%%%%%%%%%%%

\newpage
\small
\bibliography{../BFG}
\bibliographystyle{amsalpha}               

\goodbreak
\bigskip\bigskip\noindent
{\small 
Nils Berglund \\ 
{\sc CPT--CNRS Luminy} \\
Case 907, 13288~Marseille Cedex 9, France \\
{\it and} \\
{\sc PHYMAT, Universit\'e du Sud Toulon--Var} \\
{\it Present address:} \\
{\sc MAPMO--CNRS, Universit\'e d'Orl\'eans} \\
B\^atiment de Math\'ematiques, Rue de Chartres \\
B.P. 6759, 45067 Orl\'eans Cedex 2, France \\
{\it E-mail address: }{\tt berglund@cpt.univ-mrs.fr}

\bigskip\noindent
Bastien Fernandez \\ 
{\sc CPT--CNRS Luminy} \\
Case 907, 13288~Marseille Cedex 9, France \\
{\it E-mail address: }{\tt fernandez@cpt.univ-mrs.fr}

\bigskip\noindent
Barbara Gentz \\ 
{\sc Weierstra\ss\ Institute for Applied Analysis and Stochastics} \\
Mohrenstra{\ss}e~39, 10117~Berlin, Germany \\
{\it Present address:}\\
{\sc Faculty of Mathematics, University of Bielefeld} \\
P.O. Box 10 01 31, 33501~Bielefeld, Germany \\
{\it E-mail address: }{\tt gentz@math.uni-bielefeld.de}

}

%%%%%%%%%%%%%%%%%%%%%%%%%%%%%%%%%%%%%%%%%%%%%%%%%%%%%%%%%%%%%%%%%%%%%%%%%%%

\end{document}